# DEPTH WEIGHTED SCATTER ESTIMATORS

## By Yijun Zuo[1] and Hengjian Cui[2]

### *Michigan State University and Beijing Normal University*

General depth weighted scatter estimators are introduced and investigated. For general depth functions, we find out that these affine equivariant scatter estimators are Fisher consistent and unbiased for a wide range of multivariate distributions, and show that the sample scatter estimators are strong and $\sqrt{n}$-consistent and asymptotically normal, and the influence functions of the estimators exist and are bounded in general. We then concentrate on a specific case of the general depth weighted scatter estimators, the projection depth weighted scatter estimators, which include as a special case the well-known Stahel–Donoho scatter estimator whose limiting distribution has long been open until this paper. Large sample behavior, including consistency and asymptotic normality, and efficiency and finite sample behavior, including breakdown point and relative efficiency of the sample projection depth weighted scatter estimators, are thoroughly investigated. The influence function and the maximum bias of the projection depth weighted scatter estimators are derived and examined. Unlike typical high-breakdown competitors, the projection depth weighted scatter estimators can integrate high breakdown point and high efficiency while enjoying a bounded-influence function and a moderate maximum bias curve. Comparisons with leading estimators on asymptotic relative efficiency and gross error sensitivity reveal that the projection depth weighted scatter estimators behave very well overall and, consequently, represent very favorable choices of affine equivariant multivariate scatter estimators.

**1. Introduction.** The sample mean vector and sample covariance matrix have been the standard estimators of location and scatter in multivariate

Received May 2002; revised January 2004.

[1]Supported in part by NSF Grants DMS-00-71976 and DMS-01-34628.

[2]Supported in part by the NSFC and EYTP of China.

statistics. They are affine equivariant and highly efficient for normal population models. They, however, are notorious for being sensitive to unusual observations and susceptible to small perturbations in data. $M$-estimators [Maronna (1976)] are the early robust alternatives which have reasonably good efficiencies while being resistant to small perturbations in the data. Like their predecessors, the $M$-estimators unfortunately are not globally robust in the sense that they have relatively low breakdown points in high dimensions. The Stahel–Donoho (S–D) estimator [Stahel (1981) and Donoho (1982)] is the first affine equivariant estimator of multivariate location and scatter which attains a very high breakdown point. The estimator has stimulated extensive research in seeking affine equivariant location and scatter estimators which possess high breakdown points. Though $\sqrt{n}$-consistent [Maronna and Yohai (1995)], the limiting distribution of the S–D estimator remained unknown until very recently. This drawback has severely hampered the estimator from becoming more prevalent and useful in practical inference. The limiting distribution of the S–D (and general depth weighted) *location* estimator(s) has recently been discovered by Zuo, Cui and He (2004). Establishing the limiting distribution (and studying other properties) of general depth weighted and (particularly) the S–D *scatter* estimators is one goal of this paper.

In addition to the S–D estimator, affine equivariant estimators of multivariate location and scatter with high breakdown points include the minimum volume ellipsoid (MVE) and the minimum covariance determinant (MCD) estimators [Rousseeuw (1985)] and $S$-estimators [Davies (1987) and Lopuhaä (1989)]. A drawback to many classical high breakdown point estimators though is the lack of good efficiency at uncontaminated normal models. Estimators which can combine good global robustness (high breakdown point and moderate maximum bias curve) and local robustness (bounded influence function and high efficiency) are always desirable. Proposing (and investigating) a class of such estimators is another goal of this paper.

Breakdown point serves as a measure of global robustness, while the influence function captures the local robustness of estimators. In between the two extremes comes the maximum bias curve. A discussion of the maximum bias curve of scatter estimators at population models (with unknown location), seemingly very natural and desirable, has not yet been seen in the literature, perhaps partially because of the complication and difficulty to derive it. Providing an account of the maximum bias of projection depth weighted scatter estimators is the third goal of this paper.

To these ends, general depth weighted estimators are introduced and studied. The S–D estimator is just a special case of these general estimators. The paper investigates the asymptotics of the general depth weighted scatter estimators. Sufficient conditions for the asymptotic normality and the existence of influence functions of the general estimators are presented. They are satisfied by common depth functions including Tukey halfspace [Tukey (1975)] and Liu simplicial [Liu (1990)] depth. The paper then specializes to the projection depth weighted scatter estimators and examines their large and finite sample behavior. The asymptotic normality of the S–D scatter estimator follows as a special case. The influence function (together with the asymptotic relative efficiency) of the projection depth weighted scatter estimators is compared to those of some leading estimators. To fulfill the third goal of the paper, the maximum bias (under point-mass contamination) of the projection depth weighted scatter estimators for elliptical symmetric models is derived.

Findings in the paper reveal that the S–D and the projection depth weighted scatter estimators possess good robustness properties locally (high efficiency and bounded influence function) and globally (high breakdown point and moderate maximum bias) and behave very well overall compared with the leading competitors and, thus, represent favorable choices of scatter estimators.

The empirical process theory approach in the paper is useful for other depth applications. The treatment of the maximum bias of scatter estimators here sets a precedent for similar problems.

The rest of the paper is organized as follows. Section 2 introduces general depth weighted scatter estimators and investigates their Fisher consistency, asymptotics and influence functions. Section 3 is devoted to a specific case of the general depth weighted scatter estimators, the projection depth weighted scatter estimators. Here, sufficient conditions introduced in Section 2 for asymptotics and influence functions are verified and the corresponding general results are also concretized. Furthermore, the asymptotic relative efficiency, the influence function and the gross error sensitivity of the estimators are derived and compared with those of leading estimators. The maximum bias curve (under point-mass contamination) of the estimators is also derived and examined. Finally, the finite sample behavior of the estimators, including breakdown point and relative efficiency, is investigated. Simulation results with contaminated and uncontaminated data confirm the validity of the asymptotic properties at finite samples. The paper ends in Section 4 with some concluding remarks. Selected (sketches of) proofs and auxiliary lemmas are saved for the Appendix.

## 2. General depth weighted scatter estimators.

Depth functions can be employed to extend the univariate $L$-functionals ($L$-statistics) to the multivariate setting [Liu (1990) and Liu, Parelius and Singh (1999)]. For example, one can define a depth-weighted mean based on a given depth function $D(x, F)$ as follows [Zuo, Cui and He (2004)]:

$$(1) \qquad L(F) = \int x w_1(D(x, F)) \, dF(x) \Big/ \int w_1(D(x, F)) \, dF(x),$$

where $w_1(\cdot)$ is a suitable weight function [$w_1$ and $D$ are suppressed in $L(\cdot)$ for simplicity]. Subsequently, a depth-weighted scatter estimator based on $D(x, F)$ can be defined as

$$
\begin{aligned}
(2) \qquad S(F) = \int (x - L(F))(x - L(F))' \\
\times w_2(D(x, F)) \, dF(x) \Big/ \int w_2(D(x, F)) \, dF(x),
\end{aligned}
$$

where $w_2(\cdot)$ is a suitable weight function that can be different from $w_1(\cdot)$. $L(\cdot)$ and $S(\cdot)$ include multivariate versions of trimmed means and covariance matrices. The latter are excluded in later discussion though for technical convenience. To ensure well-defined $L(F)$ and $S(F)$, we require

$$
\begin{aligned}
(3) \qquad \int w_i(D(x, F)) \, dF(x) > 0, \\
\int \|x\|^i w_i(D(x, F)) \, dF(x) < \infty, \qquad i = 1, 2,
\end{aligned}
$$

where $\|\cdot\|$ stands for the Euclidean norm. The first part of (3) holds automatically for typical weight and depth functions and the second part becomes trivial if $E\|X\|^2 < \infty$ or if $w_i$, $i = 1, 2$, vanishes outside some bounded set. Replacing $F$ with its empirical version $F_n$, we obtain $L(F_n)$ and $S(F_n)$ as empirical versions of $L(F)$ and $S(F)$, respectively. $L(\cdot)$ and $S(\cdot)$ distinguish themselves from other leading estimators such as MVE- and MCD-, $S$-, $M$- and CM-estimators in the sense that $L(\cdot)$ is defined independently of $S(\cdot)$. They are also different from the ones in Lopuhaä [1999] since no prior location and scatter estimators are needed to define themselves. With the projection depth function $PD(\cdot, \cdot)$ (see Section 3), $L(\cdot)$ and $S(\cdot)$ include as special cases the well-known Stahel–Donoho location and scatter estimators, respectively.

In addition to $PD(\cdot, \cdot)$, common choices of $D(\cdot, \cdot)$ include the Tukey [1975] halfspace depth function, $HD(x, F) = \inf\{P(H) : H \text{ a closed halfspace}, x \in H\}$, and the Liu [1990] simplicial depth function, $SD(x, F) = P(x \in S[X_1, \ldots, X_{d+1}])$, where $X_1, \ldots, X_{d+1}$ is a random sample from $F$ and $S[x_1, \ldots, x_{d+1}]$ denotes the $d$-dimensional simplex with vertices $x_1, \ldots, x_{d+1}$. Weighted or trimmed means based on the latter two depth functions were considered in Liu [1990], Dümbgen [1992] and Massé [2004]. For all these depth functions, $L(\cdot)$ and $S(\cdot)$ are *affine equivariant*, that is, $L(F_{AX+b}) = AL(F) + b$, and $S(F_{AX+b}) = AS(F)A'$ for any $d \times d$ nonsingular matrix $A$ and vector $b \in \mathbb{R}^d$. In fact, this is true for any *affine invariant* $D(\cdot, \cdot)$ [i.e., $D(Ax + b, F_{AX+b}) = D(x, F)$]. With such $D(\cdot, \cdot)$ and for $F$ centrally symmetric about $\theta \in \mathbb{R}^d$ [i.e., $F_{X-\theta}(\cdot) = F_{\theta-X}(\cdot)$], $L(F)$ is *Fisher consistent* [$L(F) = \theta$] and $L(F_n)$ is *unbiased* for $\theta$ if $EX < \infty$ [Zuo, Cui and He [2004]]. This turns out to be true also for $S(F)$ and $S(F_n)$. That is, for a broad class of symmetric distributions $F$ (including as special cases elliptically symmetric $F$) with $E\|X\|^2 < +\infty$, $S(F) = \kappa\, \mathrm{Cov}(X)$ and $E(S(F_n)) = \kappa_n\, \mathrm{Cov}(X)$, for some positive constants $\kappa$ and $\kappa_n$ (with $\kappa_n \to \kappa$ as $n \to \infty$).

$L(F)$ and $L(F_n)$ have been studied in Zuo, Cui and He (2004) and Zuo, Cui and Young (2004) with respect to robustness and large and finite sample behavior. This current paper focuses on $S(F)$ and $S(F_n)$. Throughout the paper, we assume that $0 \le D(x, F) \le 1$ and $D(\cdot, \cdot)$ is continuous in $x$ and translation invariant, that is, $D(x + b, F_{X+b}) = D(x, F)$ for the given $F$ and for any $b \in \mathbb{R}^d$.

*$\sqrt{n}$-consistency and asymptotic normality.* Define

$$H_n(\cdot) = \sqrt{n}(D(\cdot, F_n) - D(\cdot, F)), \qquad \|H_n\|_\infty = \sup_{x \in \mathbb{R}^d} |H_n(x)|.$$

For a given $F$, denote $D_r = \{x : D(x, F) \ge r\}$ for $0 \le r \le 1$. Let $w_i^{(1)}$ be the derivative of $w_i$ for $i = 1, 2$. A function $g(\cdot)$ on $[a, b]$ is said to be Lipschitz continuous if there is some $C > 0$ such that $|g(s) - g(t)| \le C|s - t|$ for any $s, t \in [a, b]$. For $0 \le r_0 \le 1$, define the conditions:

(A1) $\|H_n\|_\infty = O_p(1)$ and $\sup_{x \in D_{r_0}} \|x\|\|H_n(x)| = O_p(1)$.

(A2) $w_i(r)$, $i = 1, 2$, is continuously differentiable on $[0, 1]$ and $0$ on $[0, \alpha r_0]$ for some $\alpha > 1$, $w_2^{(1)}(r)$ is Lipschitz continuous on $[0, 1]$, $w_2^{(1)}(0) = 0$, and $\int_{D_{r_0}} \|x\| |w_2^{(1)}(D(x, F))|\, dF(x) < \infty$.

In light of Vapnik–Červonenkis classes and the CLT for empirical processes [Pollard (1984) and van der Vaart and Wellner (1996)], it is seen that the first part of (A1) holds for common $D(\cdot, \cdot)$ such as $HD(\cdot, \cdot)$ and $SD(\cdot, \cdot)$. The first part of (A2) holds automatically for smooth $w_i$ such as

$$(4) \quad \begin{aligned} w_i(r) &= ((\exp(-K(1 - (r/C)^{2i})^{2i}) - \exp(-K))/(1 - \exp(-K)))I(r < C) \\ &\quad + I(r \ge C), \end{aligned}$$

with parameters $0 < C < 1$ and $K > 0$ and indicator function $I(\cdot)$ (here $r_0 = 0$), $i = 1, 2$, which will be used later. Note that (A2) excludes the trimmed means and covariance matrices with indicator functions as $w_i$. This, however, allows us to impose fewer and less severe conditions on $F$ and $D(\cdot, \cdot)$. The second part of (A1) or (A2) holds with any $r_0 > 0$ for common depth functions, in virtue of their "vanishing at infinity" property [Liu (1990) and Zuo and Serfling (2000a, b)], that is, $\lim_{\|x\| \to \infty} D(x, F) = 0$. In Section 3 we show that (A1) and (A2) hold for $PD(\cdot, \cdot)$ with $r_0 = 0$.

THEOREM 2.1. *Under* (A1) *and* (A2), $S(F_n) - S(F) = O_p(1/\sqrt{n})$.

The (strong) consistency of $S(F_n)$ can be established similarly based on corresponding conditions. Hereafter, we omit the (strong) consistency discussion. To establish the asymptotic normality of $S(F_n)$, we need the following conditions. Denote $\nu_n(\cdot) = \sqrt{n}(F_n(\cdot) - F(\cdot))$.

(A3) $\int_{D_{r_0}} \|x\|^{2i} (w_i(D(x,F)))^2 \, dF(x) < \infty$, $\int_{D_{r_0}} \|x\|^i |w_i^{(1)}(D(x,F))| \, dF(x) < \infty$, $i = 1, 2$.

(A4) $H_n(x) = \int h(x,y) \, d\nu_n(y) + o_p(1)$ uniformly on $S_n \subset D_{r_0}$, $P\{D_{r_0} - S_n\} = o(1)$, for some $h$ and $\int (\int \|y\|^i |w_i^{(1)}(D(y,F))| h(y,x)| \, dF(y))^2 \, dF(x) < \infty$, $i = 1, 2$, and $\{h(x,\cdot) : x \in S_n\}$ is a Donsker class.

Note that with a positive $r_0$, (A3) holds automatically for depth functions vanishing at infinity. (A4) holds for *HD* and *SD* with any positive $r_0$ [Dümbgen ([1992](#)) and Massé (2004)] and other depth functions. For details on a Donsker class of functions, see van der Vaart and Wellner ([1996](#)). In Section 3 we show that (A3)–(A4) hold for *PD* with $r_0 = 0$ and smooth $w_i$ [such as those in ([4](#))], $i = 1, 2$.

Let vec$(\cdot)$ be the operator which stacks the columns of a $p \times q$ matrix $M = (m_{ij})$ on the top of each other, that is, vec$(M) = (m_{11}, \ldots, m_{p1}, \ldots, m_{1q}, \ldots, m_{pq})'$. Let $M_1 \otimes M_2$ be the Kronecker product of matrices $M_1$ and $M_2$. Let $k_s(\cdot, F) = (\cdot - L_1(F))(\cdot - L_1(F))' - S(F)$. Define for $i = 1, 2$,

$$L_i(F) = \frac{\int x w_i(D(x,F)) \, dF(x)}{\int w_i(D(x,F)) \, dF(x)}, \tag{5}$$

$$\begin{aligned}
K_i(x, F) = \Bigg\{ & \int (y - L_i(F)) w_i^{(1)}(D(y,F)) h(y,x) \, dF(y) \\
& \qquad\qquad + (x - L_i(F)) w_i(D(x,F)) \Bigg\} \\
& \times \left\{ \int w_i(D(x,F)) \, dF(x) \right\}^{-1}
\end{aligned} \tag{6}$$

and

$$\begin{aligned}
& K_s(x, F) \\
& = \frac{\int k_s(y,F) w_2^{(1)}(D(y,F)) h(y,x) \, dF(y) + k_s(x,F) w_2(D(x,F))}{\int w_2(D(x,F)) \, dF(x)}.
\end{aligned} \tag{7}$$

THEOREM 2.2. *Under* (A1)–(A4), *we have*

$$S(F_n) - S(F) = \frac{1}{n} \sum_{i=1}^n (K(X_i) - E(K(X_i))) + o_p\left(\frac{1}{n}\right),$$

*where* $K(\cdot) = K_s(\cdot, F) - K_1(\cdot, F)(L_2(F) - L(F))' - (L_2(F) - L(F))(K_1(\cdot, F))'$. *Hence,*

$$\sqrt{n}(\text{vec}(S(F_n)) - \text{vec}(S(F))) \xrightarrow{d} N_{d^2}(\mathbf{0}, V),$$

*where* $V$ *is the covariance matrix of* vec$(K(X))$.

The main ideas and the outline of the proof are as follows. The key problem is to approximate

$$I_{in} = \sqrt{n}\Big( \int h_i(x) w_i(D(x, F_n)) \, dF_n(x)$$
$$- \int h_i(x) w_i(D(x, F)) \, dF(x) \Big), \qquad i = 1, 2,$$

where $h_1(x) = x - L(F)$ or $1$ and $h_2(x) = k_s(x, F)$ or $1$. The difficulty lies in the first integrand—it depends on $F_n$. By differentiability of $w_i$, there is $\theta_{in}(x)$ between $D(x, F)$ and $D(x, F_n)$ such that

$$I_{in} = \int h_i(x) w_i(D(x, F)) \, d\nu_n(x) + \int h_i(x) w_i^{(1)}(\theta_{in}(x)) H_n(x) \, dF_n(x).$$

The CLT takes care of the first term on the right-hand side. Call the second term $I_{in}^2$. Then by (A1) and (A2),

$$I_{in}^2 = \int h_i(x) w_i^{(1)}(D(x, F)) H_n(x) \, dF_n(x) + o_p(1).$$

Now by virtue of (A3) and (A4) (and, consequently, asymptotic tightness of $H_n$) and Fubini's theorem,

$$I_{in}^2 = \int \Big( \int h_i(x) w_i^{(1)}(D(x, F)) h(x, y) \, dF(x) \Big) \, d\nu_n(y) + o_p(1).$$

The desired results in Theorem 2.2 follow from the above arguments. See the Appendix for details.

*Influence function.* Now we study the influence function of $S(\cdot)$. For a given distribution $F$ in $\mathbb{R}^d$ and an $\varepsilon > 0$, the version of $F$ contaminated by an $\varepsilon$ amount of an arbitrary distribution $G$ in $\mathbb{R}^d$ is denoted by $F(\varepsilon, G) = (1 - \varepsilon)F + \varepsilon G$. The *influence function* of a functional $T$ at a given point $x \in \mathbb{R}^d$ for a given $F$ is defined as [Hampel, Ronchetti, Rousseeuw and Stahel (1986)]

$$IF(x; T, F) = \lim_{\varepsilon \to 0^+} (T(F(\varepsilon, \delta_x)) - T(F))/\varepsilon,$$

where $\delta_x$ is the point-mass probability measure at $x \in \mathbb{R}^d$. $IF(x; T, F)$ describes the relative effect (influence) on $T$ of an infinitesimal point-mass contamination at $x$, and measures the local robustness of $T$. An estimator with a bounded influence function (with respect to a given norm) is therefore robust (locally, as well as globally) and very desirable. Define for any $y \in \mathbb{R}^d$,

$$H_\varepsilon(x, y) = (D(x, F(\varepsilon, \delta_y)) - D(x, F))/\varepsilon, \qquad \|H_\varepsilon(y)\|_\infty = \sup_{x \in \mathbb{R}^d} |H_\varepsilon(x, y)|.$$

If the limit of $H_\varepsilon(x, y)$ exists as $\varepsilon \to 0^+$, then it is $IF(y; D(x, F), F)$. In the following, we assume that $IF(y; D(x, F), F)$ exists. The latter is true for the halfspace [Romanazzi (2001)], the projection [Zuo, Cui and Young (2004)], the weighted $L^p$ [Zuo (2004)] and Mahalanobis depth (MD) functions. To establish the influence function of $S(\cdot)$, we need the following condition, a counterpart of (A1). Denote by $O_y(1)$ a quantity which may depend on $y$ but is bounded as $\varepsilon \to 0$.

(A1′)  $\|H_\varepsilon(y)\|_\infty = O_y(1)$ and $\sup_{x \in D_{r_0}} \|x\| |H_\varepsilon(x, y)| = O_y(1)$.

Condition (A1′) holds for $HD$ and weighted $L^p$ depth with a positive $r_0$ and for $PD$ and $MD$ with $r_0 = 0$. Replace $h(y, x)$ in (6) and (7) by $IF(x; D(y, F), F)$ and call the resulting functions $\widetilde{K}_i(x, F)$, $i = 1, 2$, and $\widetilde{K}_s(x, F)$, respectively. We have the following:

THEOREM 2.3.  *Under* (A1′) *and* (A2),

$$IF(y; S, F) = \widetilde{K}_s(y, F) - \widetilde{K}_1(y, F)(L_2(F) - L(F))'$$
$$- (L_2(F) - L(F))(\widetilde{K}_1(y, F))'.$$

For smooth $w_i$, $i = 1, 2$, the *gross error sensitivity* of $S$: $\gamma^*(S, F) = \sup_{y \in \mathbb{R}^d} \| IF(y; S, F) \|$, where (and hereafter) "$\| \cdot \|$" stands for a selected matrix norm, is bounded if $r_0 > 0$. If $r_0 = 0$, it is also bounded if $\sup_{y \in \mathbb{R}^d} \| y^i w_i(D(y, F)) \| < \infty$, $i = 1, 2$. The latter is true for $PD$ and $MD$ and suitable $w_i$, $i = 1, 2$ [such as those in (4)].

Note that the set $D_{r_0}$ in this section could be replaced by any bounded set containing $D_{r_0}$ or the whole space $\mathbb{R}^d$, depending on the application. The latter case corresponds to $r_0 = 0$. When $r_0 > 0$, by (A2), $w_i(r) = 0$, $i = 1, 2$, for $r$ in a neighborhood of 0, corresponding to a depth trimmed (and weighted) $L(F)$ and $S(F)$ and a bounded $D_{r_0}$ for any $D(\cdot, \cdot)$ vanishing at infinity.

This section provides a general mechanism for establishing the asymptotics and the influence function of general depth weighted scatter estimators. Some of the sufficient conditions presented here might be slightly weakened in some minor aspects (e.g., for $w_1$ Lipschitz continuity suffices). Also note that results in Theorems 2.2 and 2.3 become much simpler if $w_1 = w_2$ or if $F$ is centrally symmetric since $L_2(F) = L(F)$ in these cases.

**3. Projection depth weighted and Stahel–Donoho scatter estimators.** This section is specialized to the specific case of the general depth weighted scatter estimators, the projection depth weighted or Stahel–Donoho scatter estimators.

Let $\mu$ and $\sigma$ be univariate location and scale functionals, respectively. The *projection depth* of a point $x \in \mathbb{R}^d$ with respect to a given distribution $F$ of a random vector $X \in \mathbb{R}^d$, $PD(x, F)$, is defined as [Zuo and Serfling ([2000](#)a) and Zuo ([2003](#))]

$$PD(x, F) = 1/(1 + O(x, F)), \tag{8}$$

where the *outlyingness* $O(x, F) = \sup_{\|u\|=1} (u'x - \mu(F_u))/\sigma(F_u)$, and $F_u$ is the distribution of $u'X$. Throughout our discussions $\mu$ and $\sigma$ are assumed to exist for the univariate distributions involved. We also assume that $\mu$ and $\sigma$ are affine equivariant, that is, $\mu(F_{sY+c}) = s\mu(F_Y) + c$ and $\sigma(F_{sY+c}) = |s|\sigma(F_Y)$, respectively, for any scalars $s$ and $c$ and random variable $Y \in \mathbb{R}$. Replacing $F$ with its empirical version $F_n$ based on a random sample $X_1, \ldots, X_n$, an empirical version $PD(x, F_n)$ is obtained. With $\mu$ and $\sigma$ being the median (Med) and the median absolute deviation (MAD), respectively, Liu ([1992](#)) first suggested the use of $PD(x, F_n)$ as a depth function. For motivation, examples and related discussion of ([8](#)), see Zuo ([2003](#)).

To establish the asymptotics and influence function of the projection depth weighted scatter estimators, some conditions on $\mu$ and $\sigma$ are needed. Denote by $F_{nu}$ the empirical distribution function of $\{u'X_i, i = 1, \ldots, n\}$ for any unit vector $u \in \mathbb{R}^d$.

(B1) $\sup_{\|u\|=1} |\mu(F_u)| < \infty$, $\sup_{\|u\|=1} \sigma(F_u) < \infty$ and $\inf_{\|u\|=1} \sigma(F_u) > 0$.
(B2) $\sup_{\|u\|=1} |\mu(F_{nu}) - \mu(F_u)| = O_p(1/\sqrt{n})$, $\sup_{\|u\|=1} |\sigma(F_{nu}) - \sigma(F_u)| = O_p(1/\sqrt{n})$.

Conditions (B1) and (B2) hold for common choices of $(\mu, \sigma)$ and a wide range of distributions; see Remark 2.4 of Zuo ([2003](#)) for a detailed discussion [also see Zuo, Cui and He ([2004](#))].

### 3.1. *Large sample behavior and influence function.*

#### 3.1.1. *General distributions.*

$\sqrt{n}$-*consistency and asymptotic normality.* Denote by $PWS(\cdot)$ a $PD$ weighted scatter estimator. To establish the $\sqrt{n}$-consistency of $PWS(F_n)$, we need the following lemma [Zuo ([2003](#))]:

LEMMA 3.1. *Under* (B1) *and* (B2), $\sup_{x \in \mathbb{R}^d} (1 + \|x\|)|PD(x, F_n) - PD(x, F)| = O_p(1/\sqrt{n})$.

By the lemma, (A1) holds for $PD$ with $r_0 = 0$ under (B1) and (B2). For smooth $w_i$, $i = 1, 2$, (A2) also holds since $\sup_{x \in \mathbb{R}^d} \|x\| PD(x, F) < \infty$ under (B1) [see the proof of Theorem 2.3 of Zuo ([2003](#))] and $\int \|x\| w_2^{(1)}(PD(x, F)) \, dF(x) \leq C \int \|x\| PD(x, F) \, dF(x) < \infty$. These and Theorem [2.1](#) lead to the next theorem.

THEOREM 3.1. *Assume that $w_1^{(1)}(r)$ is continuous and $w_2^{(1)}(r)$ is Lipschitz continuous on $[0, 1]$, $w_i^{(1)}(r) = O(r^i)$ for small $r \geq 0$, and $\int w_i(PD(x, F)) \, dF(x) > 0$, $i = 1, 2$. Then under (B1) and (B2), $PWS(F_n) - PWS(F) = O_p(1/\sqrt{n})$.*

Maronna and Yohai (1995) showed the $\sqrt{n}$-consistency of the S–D scatter estimator, a special case of $PWS(F_n)$ (and with $w_1 = w_2$). In Theorem 3.1 $w_i^{(1)}(r) = O(r^i)$ for small $r \geq 0$ can be relaxed to $w_i(0) = 0$ and $w_2^{(1)}(0) = 0$, $i = 1, 2$. Note that $w_i$ in (4) can serve as $w_i$ in Theorem 3.1.

For smooth $w_i$, $i = 1, 2$, in Theorem 3.1, it is readily seen that (A3) holds with $r_0 = 0$ under (B1). To establish the asymptotic normality of $PWS(F_n)$, we need to verify (A4). For any $x$ let $u(x)$ be the set of unit vectors $u$ satisfying $O(x, F) = (u'x - \mu(F_u))/\sigma(F_u)$. If $u(x)$ is a singleton, we also use $u(x)$ as the unique direction. If $X$ is a continuous random variable, nonuniqueness of $u(x)$ may occur at finitely many points. Define the following conditions:

(C1)  $\mu(F_u)$ and $\sigma(F_u)$ are continuous in $u$, $\sigma(F_u) > 0$, and $u(x)$ is a singleton except for points $x \in A \subset \mathbb{R}^d$ with $P(A) = 0$.

(C2)  The asymptotic representations $\mu(F_{nu}) - \mu(F_u) = \frac{1}{n} \sum_{i=1}^n f_1(X_i, u) + o_p(1/\sqrt{n})$ and $\sigma(F_{nu}) - \sigma(F_u) = \frac{1}{n} \sum_{i=1}^n f_2(X_i, u) + o_p(1/\sqrt{n})$ hold uniformly for $u$, the graph set of $\{f_j(X, u) : \|u\| = 1\}$ forms a polynomial set class with $E(f_j(X, u)) = 0$ for any $\|u\| = 1$,

$$E\left[ \sup_{\|u\|=1} f_j^2(X, u) \right] < +\infty$$

and

$$E\left[ \sup_{|u_1 - u_2| \leq \delta} |f_j(X, u_1) - f_j(X, u_2)|^2 \right] \to 0 \qquad \text{as } \delta \to 0, j = 1, 2.$$

For details on polynomial set classes, see Pollard (1984). (C1) and (C2) hold for general $M$-estimators of location and scale and a wide range of distributions; see Zuo, Cui and He (2004) for further discussion. Under these conditions we obtain the following [Zuo, Cui and He (2004)].

LEMMA 3.2. *Under conditions (C1) and (C2), there exists a sequence of sets $S_n \subset \mathbb{R}^d$ such that $1 - P\{S_n\} = o(1)$ and $H_n(x) = \int h(x, y) \, d\nu_n(y) + o_p(1)$ uniformly over $S_n$ with*

(9)  $h(x, y) = (O(x, F) f_2(y, u(x)) + f_1(y, u(x))) / (\sigma(F_{u(x)})(1 + O(x, F))^2).$

Hence, for smooth $w_i$, $i = 1, 2$, in Theorem 3.1, (A4) holds for $PD$ under (B1) and (C1) and (C2) with $r_0 = 0$ [see Section 2.10.2 of van der Vaart and Wellner (1996) for the verification of a Donsker class]. In light of Theorem 2.2 for general depth weighted scatter estimators, we have the following:

THEOREM 3.2. *For $w_i$, $i = 1, 2$, in Theorem* 3.1 *and under* (B1) *and* (B2) *and* (C1) *and* (C2),

$$PWS(F_n) - PWS(F) = \frac{1}{n} \sum_{i=1}^{n} K(X_i) + o_p\left(\frac{1}{n}\right),$$

*where* $K(x) = K_s(x, F) - K_1(x, F)(L_2(F) - L(F))' - (L_2(F) - L(F)) \times (K_1(x, F))'$. *Hence*

$$\sqrt{n}(\text{vec}(PWS(F_n)) - \text{vec}(PWS(F))) \xrightarrow{d} N(\mathbf{0}, V),$$

*where* $V$ *is the covariance matrix of* $\text{vec}(K(X))$.

*Influence function.* Now we derive the influence function of the projection depth weighted scatter matrices. First we need the following lemma [Zuo, Cui and Young (2004)].

LEMMA 3.3. *Assume that* (C1) *holds and the influence functions* $IF(u'y; \mu, F_u)$ *and* $IF(u'y; \sigma, F_u)$ *exist and are continuous for a given* $y \in \mathbb{R}^d$ *at* $u = u(x)$ *which is a singleton. Then*

$$
\begin{aligned}
(10) \quad & IF(y; PD(x, F), F) \\
& = \frac{O(x, F)IF((u(x))'y; \sigma, F_{u(x)}) + IF((u(x))'y; \mu, F_{u(x)})}{\sigma(F_{u(x)})(1 + O(x, F))^2}.
\end{aligned}
$$

Condition (B1) holds automatically under the conditions of this lemma and, consequently, it can be shown that (A1′) holds with $r_0 = 0$. By Theorem 2.3 we have the next theorem.

THEOREM 3.3. *Under the conditions of Lemma* 3.3 *and for smooth* $w_i$, *$i = 1, 2$, in Theorem* 3.1,

$$
\begin{aligned}
IF(y; PWS, F) = {}& \widetilde{K}_s(y, F) - \widetilde{K}_1(y, F)(L_2(F) - L(F))' \\
& - (L_2(F) - L(F))(\widetilde{K}_1(y, F))'.
\end{aligned}
$$

The influence function $IF(y; PWS, F)$ in Theorem 3.3 can be shown (details skipped) to be uniformly bounded in $y \in \mathbb{R}^d$ (with respect to a matrix norm). Thus, $\gamma^*(PWS, F) < \infty$.

3.1.2. *Elliptically symmetric distributions.* Now we focus on elliptically symmetric $F$ and $(\mu, \sigma) = (\text{Med}, \text{MAD})$. $X \sim F_{\theta, \Sigma}$ is *elliptically symmetric* about $\theta$ with a positive definite matrix $\Sigma$ *associated* if for any unit vector $u$, $u'(X - \theta) \overset{d}{=} \sqrt{u'\Sigma u} Y$ with $Y \overset{d}{=} -Y$, where "$\overset{d}{=}$" stands for "equal in distribution." First we have this lemma:

Lemma 3.4. *Let* $\mathrm{MAD}(Y) = m_0$ *and the density* $p(y)$ *of* $Y$ *be continuous with* $p(0)p(m_0) > 0$. *Then* $u(x)$ *is a singleton except at* $x = \theta$, *and* (B1) *and* (B2) *and* (C1) *and* (C2) *hold with*

$$f_1(x, u) = \sqrt{u'\Sigma u}(\tfrac{1}{2} - I\{u'(x - \theta) \leq 0\})/p(0),$$

$$f_2(x, u) = \sqrt{u'\Sigma u}(\tfrac{1}{2} - I\{|u'(x - \theta)| \leq m_0\sqrt{u'\Sigma u}\})/2p(m_0).$$

The main part of the proof is largely based on Cui and Tian ([1994](#)) and the details are skipped. Asymptotic normality (and consistency) of $PWS(F_n)$ follows immediately from this lemma and Theorem [3.2](#). The covariance matrix $V$ in Theorem [3.2](#) can be concretized.

*Asymptotic normality.* Note that $Z = \Sigma^{-1/2}(X - \theta) \sim F_0$ is spherically symmetric about the origin and $U = (U_1, \ldots, U_d)' = Z/\|Z\|$ is uniformly distributed on the unit sphere $\{x \in \mathbb{R}^d; \|x\| = 1\}$ and is independent of $\|Z\|$ [Muirhead ([1982](#))]. Define

$$s_0(x) = 1/(1 + x/m_0),$$

$$s_i(x) = E(U_1^{2(i-1)} \operatorname{sign}(|U_1|x - m_0)), \qquad i = 1, 2,$$

$$c_0 = Ew_2(\|Z\|),$$

$$c_1 = E(\|Z\|^2 w_2(\|Z\|)))/(dc_0),$$

$$c_j = E(\|Z\|^{2j-3} s_0^2(\|Z\|) w_2^{(1)}(s_0(\|Z\|)))/(4m_0^2 p(m_0)), \qquad j = 2, 3,$$

$$t_1(x) = c_3(s_2(x) - (s_1(x) - s_2(x))/(d - 1)) + x^2 w_2(s_0(x)),$$

$$t_2(x) = c_3(s_1(x) - s_2(x))/(d - 1) - c_1 c_2 s_1(x) - c_1 w_2(s_0(x)),$$

*where* $(s_1(x) - s_2(x))/(d - 1)$ *is defined to be* 0 *when* $d = 1$.

Corollary 3.1. *Under the condition of Lemma* [3.4](#) *and for* $w_i$, $i = 1, 2$, *in Theorem* [3.1](#),

$$PWS(F_n) - PWS(F) = \frac{1}{n}\sum_{i=1}^{n} K(X_i) + o_p\left(\frac{1}{n}\right)$$

*with* $K(X) = \Sigma^{1/2}(t_1(\|Z\|)UU' + t_2(\|Z\|)I_d)\Sigma^{1/2}/c_0$ *and*

$$\sqrt{n}(\operatorname{vec}(PWS(F_n)) - \operatorname{vec}(PWS(F))) \xrightarrow{d} N(\mathbf{0}, V)$$

*with* $V = \sigma_1(I_{d^2} + K_{d,d})(\Sigma \otimes \Sigma) + \sigma_2 \operatorname{vec}(\Sigma)\operatorname{vec}(\Sigma)'$, *where* $\sigma_1 = 1/(d(d + 2)c_0^2)Et_1^2(\|Z\|)$, $\sigma_2 = \sigma_1 + \frac{2}{dc_0^2}E(t_1(\|Z\|)t_2(\|Z\|)) + \frac{1}{c_0^2}Et_2^2(\|Z\|)$, *and* $K_{d,d}$ *is a* $d^2 \times d^2$-*block matrix with* $(i, j)$-*block being equal to* $\delta_{ji}$, $\delta_{ji}$ *is a* $d \times d$-*matrix which is* 1 *at entry* $(j, i)$ *and* 0 *everywhere else*, $i, j = 1, \ldots, d$.

*Asymptotic relative efficiency.* With asymptotic normality established above, we now are in a position to study the asymptotic relative efficiency of the scatter estimator $PWS(F_n)$. We shall focus on its estimation of the "shape" of $\Sigma$, that is, its "shape component"; see Tyler (1983) and Kent and Tyler (1996) for detailed arguments. For a given shape measure $\phi$, $H(\phi; PWS, F) = \phi(\Sigma^{-1/2} PWS(F) \Sigma^{-1/2})$ measures the shape (or bias) of $PWS(F)$ with respect to $\Sigma$. It clearly is affine invariant. One example of $\phi$ is the likelihood ratio test statistic $\phi_0$ measuring the ellipticity (sphericity) of any positive definite $T$ [see Muirhead (1982), also see Maronna and Yohai (1995)],

$$\phi_0(T) = (\operatorname{trace}(T)/d)^d / \det(T).$$

For this $\phi_0$, $n \log(H(\phi_0; PWS, F_n))$ has a limiting distribution. More generally, we have the following:

THEOREM 3.4. *Assume that scatter functional $S(\cdot)$ is affine equivariant and for elliptically symmetric $F_{\theta,\Sigma}$, $S(F) = c\Sigma$ for some $c > 0$ and $\sqrt{n}(\operatorname{vec}(S(F_n)) - \operatorname{vec}(S(F))) \xrightarrow{d} N(\mathbf{0}, V)$ with $V = s_1(I_{d^2} + K_{d,d})(\Sigma \otimes \Sigma) + s_2 \operatorname{vec}(\Sigma) \operatorname{vec}(\Sigma)'$, for some $s_i > 0$, $i = 1, 2$. Then*

$$n \log(\phi_0(\Sigma^{-1/2} S(F_n) \Sigma^{-1/2})) \xrightarrow{d} \frac{s_1}{c^2} \chi^2_{(d-1)(d+2)/2} \qquad \text{as } n \to \infty.$$

The details of the proof are skipped, but the main ideas are as follows. By affine equivariance of $S(\cdot)$, assume $\Sigma = I_d$. Then we can write $S(F_n) = c(I_d + n^{-1/2} Z/c)$ with $N(\mathbf{0}, V)$ as the asymptotic distribution of $\operatorname{vec}(Z)$, where $Z = (z_{ij})$. Now expand $n \log(\phi_0(\Sigma^{-1/2} S(F_n) \Sigma^{-1/2}))$ and write

$$n \log(\phi_0(\Sigma^{-1/2} S(F_n) \Sigma^{-1/2})) = (\operatorname{trace}(Z^2) - (\operatorname{trace}(Z))^2/d)/(2c^2) + O_p(n^{-1/2})$$

$$= \tilde{z}' B \tilde{z}/c^2 + O_p(n^{-1/2}),$$

with $\tilde{z} = (z_{11}/\sqrt{2}, \ldots, z_{dd}/\sqrt{2}, z_{12}, \ldots, z_{1d}, z_{23}, \ldots, z_{(d-1)d})'$ and $B = \operatorname{diag}(I_d - \mathbf{1}\mathbf{1}'/d, I_{d(d-1)/2})$, where $\mathbf{1} = (1)_{d \times 1}$. Let $A$ be the asymptotic covariance matrix of $\tilde{z}$. Then $BAB = s_1 B$. The desired result follows since the rank of $B$ is $(d-1) \times (d+2)/2$. For related discussion see Muirhead (1982).

In light of Theorem 3.4, for $PWS(F_n)$, $s_i = \sigma_i$, $i = 1, 2$, and $c = c_1$ are given in Corollary 3.1; for the sample covariance matrix $COV(F_n)$, $c = 1$ and $s_1 = 1 + \kappa$ if $F_{\theta,\Sigma}$ has kurtosis $3\kappa$ [Tyler (1982)]. Clearly, the ratio $c_1^2(1 + \kappa)/\sigma_1$ measures the asymptotic relative efficiency (ARE) of $PWS(F_n)$ with respect to $COV(F_n)$ at the given model $F_{\theta,\Sigma}$. The same idea was employed in Tyler (1983) to compute AREs of scatter estimators. At the multivariate normal model, $\kappa = 0$, hence the ratio $c_1^2/\sigma_1$ is the ARE of $PWS(F_n)$ with respect to $COV(F_n)$.

Consider $w_i$, $i = 1, 2$, in (4). They are selected to meet the requirements in Theorem 3.1 and to down-weight exponentially less deep points to get better performance of PWS. Also, appropriate tuning of $C$ and $K$ can lead to highly efficient (and robust) PWS [see Zuo, Cui and He (2004) for related comments]. The behavior of $w_2$ is depicted in Figure 1 with $C = 0.32$ and $K = 0.2$.

Table 1 reports the AREs of $PWS(F_n)$ [with respect to $COV(F_n)$] versus the dimension $d$ and selected $C$ and $K$ at $N(\mathbf{0}, I_d)$ with $w_2$ above. Here we select $C$'s that are close to $\text{Med}(PD(X, F))$ to get better performance of PWS. It is seen that $PWS(F_n)$ possesses very high ARE for suitable $K$ and $C$, which, in fact, approaches 100% rapidly as the dimension $d$ increases. Note that the ARE of $PWS(F_n)$ here does not depend on that of the underlying projection depth weighted mean (PWM). The ARE of the latter depends on $w_1$ and behaves like that of $PWS(F_n)$ [Zuo, Cui and He (2004)].

*Influence function.* Under the condition of Lemma 3.4, it can be shown that

$$IF(u(x)'y, \text{Med}, F_{u(x)}) = \frac{\|\Sigma^{-1/2}x\|}{2p(0)\|\Sigma^{-1}x\|}\,\text{sign}(x'\Sigma^{-1}y),$$

$$IF(u(x)'y, \text{MAD}, F_{u(x)}) = \frac{\|\Sigma^{-1/2}x\|}{4p(m_0)\|\Sigma^{-1}x\|}\,\text{sign}(|x'\Sigma^{-1}y| - m_0\|\Sigma^{-1/2}x\|).$$

These functions are continuous at $u(x)$ almost surely. By Lemmas 3.4 and 3.3 we have

$$IF(x; PD(y, F), F)$$

$$= \frac{s_0^2(\|\Sigma^{-1/2}y\|)}{m_0}$$

$$\times \left(\frac{\|\Sigma^{-1/2}y\|\,\text{sign}(|y'\Sigma^{-1}x| - m_0\|\Sigma^{-1/2}y\|)}{4m_0p(m_0)} + \frac{\text{sign}(y'\Sigma^{-1}x)}{2p(0)}\right).$$

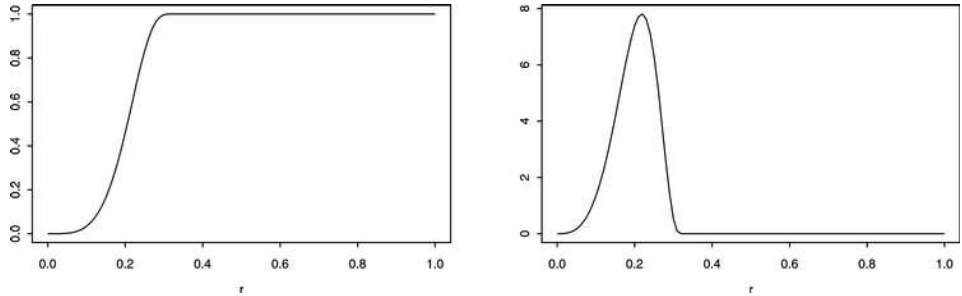

Fig. 1. *The behavior of $w_2(r)$ with $C = 0.32$ and $K = 0.2$. Left: $w_2(r)$. Right: $w_2^{(1)}(r)$.*

Table 1
*The asymptotic relative efficiency of PWS versus the dimension d*

| $d$ | $C = \frac{1}{1+\sqrt{d}/\Phi^{-1}(3/4)}$ $K = 2$ | $C = \frac{1}{1+\sqrt{d}/\Phi^{-1}(3/4)}$ $K = 3$ | $C = \frac{1}{1+\sqrt{2d}}$ $K = 2$ | $C = \frac{1}{1+\sqrt{2d}}$ $K = 3$ |
|---|---|---|---|---|
| 2  | 0.922 | 0.883 | 0.904 | 0.862 |
| 3  | 0.957 | 0.933 | 0.945 | 0.918 |
| 4  | 0.976 | 0.959 | 0.969 | 0.945 |
| 5  | 0.980 | 0.974 | 0.979 | 0.965 |
| 6  | 0.989 | 0.980 | 0.983 | 0.974 |
| 7  | 0.990 | 0.986 | 0.986 | 0.980 |
| 8  | 0.993 | 0.991 | 0.991 | 0.985 |
| 9  | 0.994 | 0.992 | 0.992 | 0.987 |
| 10 | 0.995 | 0.993 | 0.994 | 0.980 |
| 15 | 0.998 | 0.998 | 0.996 | 0.995 |
| 20 | 1.00  | 0.999 | 0.999 | 0.997 |
| 30 | 1.00  | 1.00  | 1.00  | 0.999 |

By virtue of Theorem 3.3, we have the next corollary.

COROLLARY 3.2.    *Under the condition of Lemma 3.4 and for $w_i$, $i = 1, 2$, in Theorem 3.1,*

$$IF(x; PWS, F_{0,I_d}) = (t_1(\|x\|)xx'/\|x\|^2 + t_2(\|x\|)I_d)/c_0,$$

$$IF(x; PWS, F_{\theta,\Sigma}) = \Sigma^{1/2}(IF(\Sigma^{-1/2}(x-\theta); PWS, F_{0,I_d}))\Sigma^{1/2}.$$

Figure 2 indicates $IF(x; PWS, F_{\theta,\Sigma})$ is uniformly bounded in $x \in \mathbb{R}^d$ relative to a matrix norm.

Maintaining a good balance between high efficiency and a bounded influence function is always a legitimate concern for estimators. Many existing high breakdown estimators fail to do so though. CM- [Kent and Tyler (1996)] and $\tau$- [Lopuhaä (1999)] estimators are among the few exceptions. In light of these papers, we consider a gross error sensitivity index for the shape of the scatter estimator $S$,

$$G_2(S, F) = \text{GES}(S, F)/((1 + 2/d)(1 - 1/d)^{1/2}),$$

where $\text{GES}(S, F)$ is the gross-error-sensitivity of $S(F)/\text{trace}(S(F))$, the shape component of the scatter functional $S(F)$. In our case it is seen that $G_2(PWS, F) = \sup_{r \geq 0} t_1(r)/(c_0(d + 2))$. Table 2 reports the ARE and $G_2$ of scatter estimators (along with those of the corresponding location estimators listed in parentheses; in the location case $G_2 = \gamma^*$) for $d = 2, 5$ and 10.

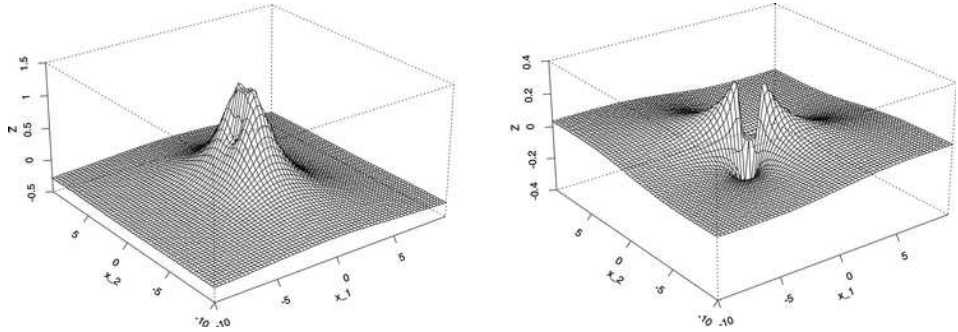

FIG. 2.    *The behavior of $IF(x; PWS, F_{0,I_2})$ with $w_2$ in (4). Left: $-(1,1)$ entry. Right: $-(1,2)$ entry.*

| $d$ | Estimator | ARE | | $G_2$ | |
|-----|-----------|-----|-----|-------|-----|
| 2 | $\tau$(CM)- | 0.8670 | (0.9057) | 1.415 | (1.861) |
|   | PWS | 0.8810 | (0.9152) | 1.318 | (1.818) |
| 5 | $\tau$(CM)- | 0.9099 | (0.9354) | 1.275 | (2.588) |
|   | PWS | 0.9180 | (0.9516) | 1.057 | (2.546) |
| 10 | $\tau$(CM)- | 0.9505 | (0.9606) | 1.224 | (3.425) |
|   | PWS | 0.9620 | (0.9734) | 0.979 | (3.421) |

Table 2 lists only the ARE and $G_2$ for $\tau$- and PWS estimators. The corresponding indices for the CM-estimators are omitted since they are almost the same as those of the $\tau$-estimators. The indices for $\tau$(CM)-estimators are obtained by optimizing $G_2$ of the corresponding location estimators based on Tukey's biweight function [Kent and Tyler (1996) and Lopuhaä (1999)]. The weight function $w_2$ in (4) is employed in our calculation for the indices of PWS [and $w_1$ in (4) for PWM] with $K = 3$ and $C = 1/(1 + \sqrt{\xi_d d})$, where $\xi_2 = 2.3$, $\xi_5 = 1.2$ and $\xi_{10} = 0.9$ for PWS [and $\xi_d = 1.2$ for PWM]. The values of $C$ here are slightly different from those in Table 1 to get (nearly) optimal ARE and $G_2$ simultaneously. Inspecting Table 2 reveals that, compared with leading competitors, the projection depth weighted scatter estimator PWS behaves very well overall.

*Maximum bias.* Define the *maximum bias* of a scatter matrix $S$ under an $\varepsilon$ amount of contamination at $F$ as $B(\varepsilon; S, F) = \sup_G \|S(F(\varepsilon, G)) - S(F)\|$, where $G$ is any distribution in $\mathbb{R}^d$. The *contamination sensitivity* of $S$ at $F$ is defined as $\gamma(S, F) = \lim_{\varepsilon \to 0+} \sup_G \|(S(F(\varepsilon, G)) - S(F))/\varepsilon\|$; see He and Simpson (1993) for a related definition for location estimators. $B(\varepsilon; S, F)$ is the maximum deviation (bias) of $S$ under an $\varepsilon$ amount of contamination at $F$, and measures mainly the global robustness of $S$. $\gamma(S, F)$ indicates the maximum relative effect on $S$ of an infinitesimal contamination at $F$, and measures the local, as well as global, robustness of $S$. The minimum amount $\varepsilon^*$ of contamination at $F$ which leads to an unbounded $B(\varepsilon; S, F)$ is called the (asymptotic) *breakdown point* (BP) of $S$ at $F$, that is, $\varepsilon^* = \min\{\varepsilon : B(\varepsilon; S, F) = \infty\}$.

In many cases, the maximum bias is attained by a point-mass distribution; see Huber ([1964]), Martin, Yohai and Zamar ([1989]), Chen and Tyler ([2002]) and Zuo, Cui and Young ([2004]). In the following, we derive the maximum bias and contamination sensitivity of the shape component of PWS under point-mass contamination. We conjecture that our results hold for general contamination. For any $0 \le \varepsilon < 1/2$ and $c \in \mathbb{R}$, define $d_1 = d_1(\varepsilon)$, $m_i(c, \varepsilon)$, $i = 1, 2$, by

$$P(Y \le d_1(\varepsilon)) = \frac{1}{2(1-\varepsilon)},$$

$$P(|Y - c| \le m_1(c, \varepsilon)) = \frac{1 - 2\varepsilon}{2(1-\varepsilon)},$$

$$P(|Y - c| \le m_2(c, \varepsilon)) = \frac{1}{2(1-\varepsilon)}$$

(assume that $d_1, m_1, m_2$ are well defined). For $x \in \mathbb{R}^d$, write $x' = (x_1, x_2')$ with $x_1 = x_{11} \in \mathbb{R}$ and $x_2 = (x_{21}, \ldots, x_{2(d-1)})' \in \mathbb{R}^{d-1}$. Likewise, partition the unit vector $u \in \mathbb{R}^d$. For any $r \ge 0$, define

$$f_1(x, r, \varepsilon) = \sup_{0 \le u_1 \le 1} \frac{\sqrt{1 - u_1^2}\, \|x_2\| + |u_1 x_1 - f_4(u_1, r, d_1)|}{f_3(u_1, r, d_1)},$$

$$f_2(r, \varepsilon) = \sup_{0 \le u_1 \le 1} \frac{|u_1 r - f_4(u_1, r, d_1)|}{f_3(u_1, r, d_1)},$$

with $f_3(u_1, r, d_1)$ being the median of $\{m_1(f_4(u_1, r, d_1), \varepsilon), |u_1 r - f_4(u_1, r, d_1)|, m_2(f_4(u_1, r, d_1), \varepsilon)\}$, $f_4(u_1, r, d_1)$ being the median of $\{-d_1, u_1 r, d_1\}$ ($\varepsilon$ is suppressed in $f_3$ and $f_4$). Define, for $i = 1, 2$,

$$\phi_i(r, \varepsilon) = (1 - \varepsilon) \int x_1 w_i \left( \frac{1}{1 + f_1(x, r, \varepsilon)} \right) dF_0(x),$$

$$\psi_i(r, \varepsilon) = (1 - \varepsilon) \int x_{i1}^2 w_2 \left( \frac{1}{1 + f_1(x, r, \varepsilon)} \right) dF_0(x),$$

$$\eta_i(r, \varepsilon) = (1 - \varepsilon) \int w_i \left( \frac{1}{1 + f_1(x, r, \varepsilon)} \right) dF_0(x),$$

$$\gamma_i(r, \varepsilon) = \varepsilon w_i \left( \frac{1}{1 + f_2(r, \varepsilon)} \right),$$

$$b_1(r, \varepsilon) = \frac{\psi_1(r, \varepsilon) - \psi_2(r, \varepsilon) + \gamma_2(r, \varepsilon)r^2}{\eta_2(r, \varepsilon) + \gamma_2(r, \varepsilon)} + \frac{(\phi_1(r, \varepsilon) + \gamma_1(r, \varepsilon)r)^2}{(\eta_1(r, \varepsilon) + \gamma_1(r, \varepsilon))^2}$$

$$- 2\frac{\phi_1(r, \varepsilon)\phi_2(r, \varepsilon) + (\phi_1(r, \varepsilon)\gamma_2(r, \varepsilon) + \phi_2(r, \varepsilon)\gamma_1(r, \varepsilon))r}{(\eta_1(r, \varepsilon) + \gamma_1(r, \varepsilon))(\eta_2(r, \varepsilon) + \gamma_2(r, \varepsilon))}$$

$$- 2\frac{\gamma_1(r, \varepsilon)\gamma_2(r, \varepsilon)r^2}{(\eta_1(r, \varepsilon) + \gamma_1(r, \varepsilon))(\eta_2(r, \varepsilon) + \gamma_2(r, \varepsilon))},$$

$$b_2(r, \varepsilon) = \psi_2(r, \varepsilon)/(\eta_2(r, \varepsilon) + \gamma_2(r, \varepsilon)) - c_1.$$

For any $y \in \mathbb{R}^d$, denote $\tilde{y} = \Sigma^{-1/2}(y - \theta)$. We have the next theorem:

**Theorem 3.5.** *Under the condition of Lemma* 3.4 *and for any $\varepsilon > 0$ and $y \in \mathbb{R}^d$,*

$$PWS(F(\varepsilon, \delta_y)) - PWS(F) = \Sigma^{1/2}(b_1(\|\tilde{y}\|, \varepsilon)\tilde{y}\tilde{y}'/\|\tilde{y}\|^2 + b_2(\|\tilde{y}\|, \varepsilon)I_d)\Sigma^{1/2}.$$

For weight functions $w_i$, $i = 1, 2$, in Theorem 3.1, it can be shown that for any $\varepsilon < 1/2$, $\text{trace}(PWS(F(\varepsilon, \delta_y)) - PWS(F))$ is uniformly bounded with respect to $y \in \mathbb{R}^d$. Hence we have the following:

**Corollary 3.3.** *Under the condition of Lemma* 3.4 *and for weight functions $w_i$, $i = 1, 2$, in Theorem* 3.1, *$\varepsilon^*(PWS, F) = 1/2$.*

Focusing again on the shape component of PWS and based on the result in Theorem 3.5, we can define in a straightforward fashion a gross error sensitivity index (GESI), a maximum bias index (MBI) and a contamination sensitivity index (CSI), respectively, as follows:

$$\text{GESI}(PWS, F) = \sup_{y \in \mathbb{R}^d} \left\|\lim_{\varepsilon \to 0} b_1(\|\tilde{y}\|, \varepsilon)\Sigma^{1/2}(\tilde{y}\tilde{y}'/\|\tilde{y}\|^2)\Sigma^{1/2}/\varepsilon\right\|,$$

$$\text{MBI}(\varepsilon; PWS, F) = \sup_{y \in \mathbb{R}^d} \|b_1(\|\tilde{y}\|, \varepsilon)\Sigma^{1/2}(\tilde{y}\tilde{y}'/\|\tilde{y}\|^2)\Sigma^{1/2}\|,$$

$$\text{CSI}(PWS, F) = \lim_{\varepsilon \to 0^+} \sup_{y \in \mathbb{R}^d} \|b_1(\|\tilde{y}\|, \varepsilon)\Sigma^{1/2}(\tilde{y}\tilde{y}'/\|\tilde{y}\|^2)\Sigma^{1/2}/\varepsilon\|.$$

In view of Corollary 3.2, it can be seen that $\text{GESI}(PWS, F) = \lambda_1 \times \sup_{r \geq 0} |t_1(r)|/c_0$, which is $\leq \text{CSI}(PWS, F)$, where $\lambda_1$ is the largest eigenvalue of $\Sigma$. Note that under point-mass contamination the only difference between CSI and GESI is the order in which the suprema and the limits are taken in their respective definitions above. This might tempt one to believe that these two sensitivity indices are the same if it is taken for granted that the order in which the supremum and the limit are taken is interchangeable. Unfortunately, this is not always the case [see, e.g., Chen and Tyler (2002)].

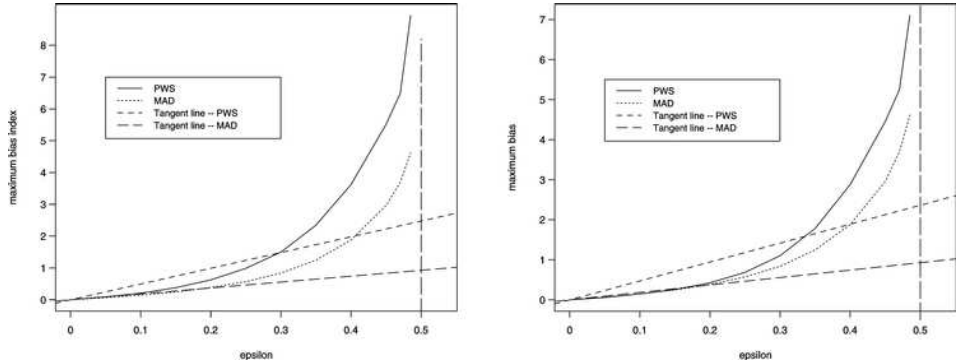

 *The behavior of the maximum bias (index) of PWS and MAD.* Left: *maximum bias indices of PWS and MAD.* Right: *maximum biases of PWS and MAD.*

In the following, we prove that for PWS, the order *is* interchangeable and $\mathrm{CSI}(PWS, F)$ is the same as $\mathrm{GESI}(PWS, F)$. The proof and the derivation of the following result, given in the Appendix, is rather technically demanding and has no precedent in the literature.

THEOREM 3.6. *Under the condition of Lemma* 3.4 *and for* $w_i$, $i = 1, 2$, *in Theorem* 3.1:

(a) $\mathrm{MBI}(\varepsilon; PWS, F) = \lambda_1 \sup_{r \geq 0} b_1(r, \varepsilon)$ *and*
(b) $\mathrm{CSI}(PWS, F) = \mathrm{GESI}(P\widetilde{W}S, F) = \lambda_1 \sup_{r \geq 0} |t_1(r)|/c_0$.

The behavior of $\mathrm{MBI}(\varepsilon; PWS, N(\mathbf{0}, I_2))$ [and $B(\varepsilon; PWS, N(\mathbf{0}, I_2))$], together with that of the (explosion) maximum bias of MAD at $N(0, 1)$ − $B(\varepsilon; \mathrm{MAD}, N(0, 1))$ (note that no separate shape and scale components correspond to MAD, a univariate scale measure), as functions of $\varepsilon$ is revealed in Figure 3. The slopes of the tangent lines at the origin represent the CSI (or $\gamma$) of PWS and MAD. From the figures we see that the maximum bias (index) of PWS is quite moderate (and slightly larger than that of the univariate scale measure MAD) and it increases very slowly as the amount of contamination $\varepsilon$ increases and jumps to infinity as $0.45 < \varepsilon \to \frac{1}{2}$, confirming that the asymptotic breakdown point of PWS is $\frac{1}{2}$.

3.2. *Finite sample behavior.* In this section the finite sample robustness and relative efficiency of $PWS(F_n)$ are investigated. Finite sample results in this section confirm the asymptotic results in the last section.

3.2.1. *Finite sample breakdown point.* Let $X^n = \{X_1, \ldots, X_n\}$ be a sample of size $n$ from $X$ in $\mathbb{R}^d$ $(d \geq 1)$. The replacement breakdown point (RBP)

[Donoho and Huber ([1983](#))] of a scatter estimator $V$ at $X^n$ is defined as

$$\text{RBP}(V, X^n) = \min\left\{\frac{m}{n} : \text{trace}(V(X^n)V(X_m^n)^{-1} + V(X^n)^{-1}V(X_m^n)) = \infty\right\},$$

where $X_m^n$ is a contaminated sample resulting from replacing $m$ points of $X^n$ with arbitrary values.

In the following discussion of the RBP of the projection depth weighted scatter estimators, $(\mu, \sigma) = (\text{Med}, \text{MAD}_k)$, where $\text{MAD}_k$ is a modified MAD which can lead to a slightly higher RBP. Similar ideas of modifying MAD to achieve higher RBP were used in Tyler ([1994](#)) and Gather and Hilker ([1997](#)). Here $\text{MAD}_k(x^n) = \text{Med}_k(\{|x_1 - \text{Med}(x^n)|, \ldots, |x_n - \text{Med}(x^n)|\})$, with $\text{Med}_k(x^n) = (x_{(\lfloor (n+k)/2 \rfloor)} + x_{(\lfloor (n+1+k)/2 \rfloor)})/2$, for $1 \le k \le n$, and $x_{(1)} \le \cdots \le x_{(n)}$ being ordered values of $x_1, \ldots, x_n$ in $\mathbb{R}^1$ (note $\text{MAD}_1 = \text{MAD}$). Denote by $PWS_n^k$ the corresponding scatter estimator.

A random sample $X^n$ is said to be *in general position* if there are no more than $d$ sample points of $X^n$ lying in any $(d-1)$-dimensional subspace. Let $\lfloor \cdot \rfloor$ be the floor function. We have the next theorem.

THEOREM 3.7. *Let $(\mu, \sigma) = (\text{Med}, \text{MAD})$ and $PD(x, F)$ be the depth function. Let $w_i(r)$ be continuous on $[0, 1]$ and positive and $\le M_i r^i$ on $(0, 1]$ for some $M_i > 0$, $i = 1, 2$. Then for $X^n$ in general position $(n > 2d)$,* $\text{RBP}(PWS_n^k, X^n) = \min\{\lfloor (n-k+2)/2 \rfloor / n, \lfloor (n+k+1-2d)/2 \rfloor / n\}.$

When $k = d$ or $d + 1$, $\text{RBP}(PWS_n^k, X^n) = \lfloor (n-d+1)/2 \rfloor / n$, the upper bound of RBP of any affine equivariant scatter estimators; see Davies ([1987](#)). The RBP of the Stahel–Donoho scatter estimator, a special case of $PWS_n^k$, has been given in Tyler ([1994](#)). Note that for the smooth $w_i$ in (A2), $w_i(r) \le M_i r^i$ holds automatically, $i = 1, 2$. The result in Theorem 3.7 holds true for any $\mu$ and $\sigma$ that share the RBPs of Med and $\text{MAD}_k$, respectively.

3.2.2. *Finite sample relative efficiency.* We generate 400 samples from the model $(1 - \varepsilon)N(\mathbf{0}, I_2) + \varepsilon\delta_{(100,0)}$ with $\varepsilon = 0\%$, 10% and 20% for sample sizes $n = 100, 200, \ldots, 1000$. An approximate algorithm with time complexity $O(n^3)$ (for $d = 2$) is utilized for the computation of the $PD_n(X_i)$, $i = 1, \ldots, n$, and the projection depth weighted scatter matrix. $(\mu, \sigma) = (\text{Med}, \text{MAD})$ and the weight functions $w_i(\cdot)$ defined in ([4](#)), with $C = 1/(1 + \sqrt{2}/\Phi^{-1}(3/4)) \approx 0.323$ and $K = 2$, are used in our simulation.

We calculate for a scatter estimator $V_n$ the mean of the likelihood ratio test (LRT) statistic $\text{LRT}(V_n) = \frac{1}{m}\sum_{j=1}^{m}\phi_0(V_j)$ with $m = 400$ and $V_j$ being the estimate for the $j$th sample. In the case with $\varepsilon = 0\%$ (no contamination), the mean of the $n$ log likelihood ratio test (LLRT) statistic with $\text{LLRT}(V_n) = \frac{1}{m}\sum_{j=1}^{m}n\log(\phi_0(V_j))$ is calculated. The finite sample relative efficiency (RE) of $V_n$ at $\varepsilon = 0\%$ is then obtained by dividing the LLRT of the sample covariance matrix by that of $V_n$ [Maronna and Yohai (1995) used the same measure for finite sample relative efficiency]. Some simulation results are listed in Table 3.

The finite sample RE of $PWS(F_n)$ related to the sample covariance matrix at $N(\mathbf{0}, I_2)$ increases from about 80% for $n = 20$ to 91% for $n = 100$ and is around 90%–93% and very stable for $n = 100, 200, \ldots, 1000$ [and is very close to its asymptotic value 92.2% (listed in Table 1)]. In the contamination cases, the results in Table 3 indicate that $PWS(F_n)$ is very robust, whereas $COV(F_n)$ is very sensitive to outliers. For the special case of $PWS_n$, the Stahel–Donoho estimator, a related simulation study was conducted by Maronna and Yohai (1995).

Though alternatives exist, the $w_2$ we select results in a very good performance of $PWS_n$ and satisfies all the requirements in the previous sections. Note that smaller $C$ can lead to a higher RE of $PWS_n$ under no contamination, while larger $C$ can lead to a better performance of $PWS_n$ under contamination. The same is true for the parameter $K$. Moderate values of $C$ and $K$ thus are recommended (and are used in our simulation); see Zuo, Cui and He (2004) for related discussion.

## 4. Concluding remarks.

General depth weighted scatter estimators are introduced and studied. The estimators possess nice properties. In a very general setting, consistency and asymptotic normality of the estimators are established and their influence functions are derived. These general results are concretized and demonstrated via the projection depth weighted scatter estimators. The latter estimators include as a special case the Stahel–Donoho estimator, the first one constructed which combines affine equivariant and high breakdown point, but has an unknown limiting distribution until this paper.

TABLE 3
*Mean of the likelihood ratio test statistic and relative efficiency*

| $n$ | PWS | COV | PWS | COV | PWS | COV | RE |
|---|---|---|---|---|---|---|---|
| | $\varepsilon = 0\%$ | | $\varepsilon = 10\%$ | | $\varepsilon = 20\%$ | | $(\varepsilon = 0\%)$ |
| 100 | 1.022 | 1.021 | 1.110 | 234.09 | 1.523 | 420.80 | 0.913 |
| 200 | 1.011 | 1.010 | 1.109 | 231.03 | 1.534 | 407.10 | 0.911 |
| 300 | 1.007 | 1.006 | 1.106 | 230.04 | 1.528 | 405.72 | 0.900 |
| 400 | 1.006 | 1.005 | 1.105 | 227.79 | 1.539 | 404.13 | 0.903 |
| 500 | 1.004 | 1.004 | 1.103 | 227.18 | 1.555 | 404.43 | 0.901 |
| 600 | 1.004 | 1.003 | 1.105 | 227.26 | 1.560 | 404.78 | 0.917 |
| 700 | 1.003 | 1.003 | 1.103 | 227.37 | 1.545 | 403.20 | 0.930 |
| 800 | 1.003 | 1.002 | 1.104 | 226.28 | 1.555 | 404.00 | 0.932 |
| 900 | 1.002 | 1.002 | 1.103 | 226.27 | 1.549 | 401.45 | 0.923 |
| 1000 | 1.002 | 1.002 | 1.102 | 226.19 | 1.543 | 401.75 | 0.926 |

Frequently high breakdown point affine equivariant estimators suffer from a low asymptotic relative efficiency and an unbounded influence function. The projection depth weight scatter estimators are proven to be exceptions. They combine the best possible breakdown point and a moderate maximum bias curve (global robustness) and a bounded influence function (local robustness) and possess, in the meantime, a very high asymptotic relative efficiency at multivariate normal models. Simulations with clean and contaminated data sets reveal that the global robustness and high efficiency properties hold at finite samples.

Finally, we comment that the $w_i$ in this paper do not include indicator functions. This allows us to treat general depth and distribution functions. To cover trimmed means (with indicator weight functions), one has to impose more conditions on these functions (but the efficiency will be lower).

## APPENDIX: SELECTED (SKETCHES OF) PROOFS AND AUXILIARY LEMMAS

PROOF OF THEOREM 2.1.   Denote by $l_1(F)$ and $l_2(F)$ the numerator and the denominator of $L(F)$, respectively, and $s_1(F)$ and $s_2(F)$ those of $S(F)$, respectively. Write

$$(11) \quad L(F_n) - L(F) = ((l_1(F_n) - l_1(F)) - L(F)(l_2(F_n) - l_2(F)))/l_2(F_n),$$

$$
\begin{aligned}
(12) \quad S(F_n) - S(F) = &\left( \int xx' w_2(D(x, F_n))\, dF_n(x) \right. \\
&\left. - \int xx' w_2(D(x, F))\, dF(x) \right) \Big/ s_2(F_n) \\
&- S_0(F)(s_2(F_n) - s_2(F))/s_2(F_n) \\
&- (L_2(F_n) - L_2(F))(L(F_n))' - L_2(F)(L(F_n) - L(F))' \\
&- (L(F))(L_2(F_n) - L_2(F))' \\
&- (L(F_n) - L(F))(L_2(F_n))' \\
&+ (L(F_n) - L(F))(L(F_n))' + L(F)(L(F_n) - L(F))',
\end{aligned}
$$

with $S_0(F) = \int xx' w_2(D(x, F))\, dF(x)/s_2(F)$. We now show that under (A1) and (A2),

$$
\begin{aligned}
(13) \quad I_{in} = &\int \|x\|^i |w_i(D(x, F_n)) - w_i(D(x, F))|\, dF_n(x) \\
&= O_p(1/\sqrt{n}), \qquad i = 1, 2.
\end{aligned}
$$

By (A2), there exists a $\theta_{in}(x)$ between $D(x, F_n)$ and $D(x, F)$ such that for $i = 1, 2$,

$$I_{in} \leq \int \|x\|^i |w_i^{(1)}(\theta_{in}(x)) - w_i^{(1)}(D(x, F))| \frac{|H_n(x)|}{\sqrt{n}} \, dF_n(x)$$

$$+ \int \|x\|^i w_i^{(1)}(D(x, F)) \frac{|H_n(x)|}{\sqrt{n}} \, dF_n(x).$$

Call the two terms in the right-hand side $I_{in}^{(1)}$ and $I_{in}^{(2)}$, respectively. Let $r_1 = \alpha r_0$. By (A1), $D(x, F) + \sup_{x \in \mathbb{R}^d} |D(x, F_n) - D(x, F)| = D(x, F) + O_p(1/\sqrt{n}) \geq \theta_{in}(x)$. This and (A2) and (A1) lead to

$$I_{in}^{(1)} = \int_{\{\theta_{in}(x) > r_1\} \cup D_{r_1}} \|x\|^i |w_i^{(1)}(\theta_{in}(x)) - w_i^{(1)}(D(x, F))| \frac{|H_n(x)|}{\sqrt{n}} \, dF_n(x)$$

$$\leq C \int_{\{D(x, F) + O_p(1/\sqrt{n}) > r_1\} \cup D_{r_1}} \|x\|^i \left( \frac{|H_n(x)|}{\sqrt{n}} \right)^i dF_n(x) = O_p\left( \left( \frac{1}{\sqrt{n}} \right)^i \right)$$

and

$$I_{in}^{(2)} = \int_{D_{r_0}} \|x\|^i w_i^{(1)}(D(x, F)) \frac{|H_n(x)|}{\sqrt{n}} \, dF_n(x) = O_p\left( \frac{1}{\sqrt{n}} \right).$$

Hence $I_{in} = O_p(1/\sqrt{n})$. Likewise we can show that

$$(14) \qquad \int w_i(D(x, F_n)) \, dF_n(x) - \int w_i(D(x, F)) \, dF(x) = O_p(1/\sqrt{n}).$$

Let $h(x) = xx'$, $x$ or 1. It follows from displays (13) and (14) and the CLT that

$$\int h(x) w_i(D(x, F_n)) \, dF_n(x) - \int h(x) w_i(D(x, F)) \, dF(x) = O_p(1/\sqrt{n}).$$

By (11), the boundedness of $L(F)$ and $l_2(F)$, and the fact that $l_2(F_n) = l_2(F) + O_p(1/\sqrt{n})$, we have $L(F_n) - L(F) = O_p(1/\sqrt{n})$. Likewise we have $L_2(F_n) - L(F) = O_p(1/\sqrt{n})$. These, (12) and the boundedness of $S_0(F)$, $s_2(F)$, $L(F)$ and $L_2(F)$ yield $S(F_n) - S(F) = O_p(1/\sqrt{n})$. $\square$

PROOF OF THEOREM 2.2. Employing the notation in the proof of Theorem 2.1, write

$$\sqrt{n} \left( \int xx' w_2(D(x, F_n)) \, dF_n(x) - \int xx' w_2(D(x, F)) \, dF(x) \right)$$

$$= \int xx' w_2^{(1)}(\theta_{2n}(x)) H_n(x) \, dF_n(x) + \int xx' w_2(D(x, F)) \, d\nu_n(x),$$

where $\theta_{2n}(x)$ is a point between $D(x, F_n)$ and $D(x, F)$. Following the proof of Theorem 2.1 and by (A1)–(A4) (and, consequently, the asymptotic tightness of $H_n$ on $S_n$), we can show that

$$\int xx' w_2^{(1)}(\theta_{2n}(x)) H_n(x) \, dF_n(x)$$
$$= \int xx' w_2^{(1)}(D(x, F)) H_n(x) \, dF_n(x) + o_p(1)$$
$$= \int xx' w_2^{(1)}(D(x, F)) H_n(x) \, dF(x) + o_p(1).$$

Therefore,

$$\sqrt{n} \left( \int xx' w_2(D(x, F_n)) \, dF_n(x) - \int xx' w_2(D(x, F)) \, dF(x) \right)$$
$$= \int xx' w_2^{(1)}(D(x, F)) H_n(x) \, dF(x) + \int xx' w_2(D(x, F)) \, d\nu_n(x) + o_p(1).$$

By (A4) and Fubini's theorem, we have

$$\sqrt{n} \left( \int xx' w_2(D(x, F_n)) \, dF_n(x) - \int xx' w_2(D(x, F)) \, dF(x) \right)$$
$$(15) \qquad = \int \left( \int yy' w_2^{(1)}(D(y, F)) h(y, x) \, dF(y) + xx' w_2(D(x, F)) \right) d\nu_n(x)$$
$$+ o_p(1).$$

Likewise, we can show that

$$\sqrt{n}(s_2(F_n) - s_2(F))$$
$$= \sqrt{n} \left( \int w_2(D(x, F_n)) \, dF_n(x) - \int w_2(D(x, F)) \, dF(x) \right)$$
$$(16) \qquad = \int \left( \int w_2^{(1)}(D(y, F)) h(y, x) \, dF(y) + w_2(D(x, F)) \right) d\nu_n(x)$$
$$+ o_p(1),$$

and for $i = 1, 2$ [see the proof of Theorem 2.1 of Zuo, Cui and He (2004)],

$$\sqrt{n}(L_i(F_n) - L_i(F))$$
$$= \left\{ \int \left( \int (y - L_i(F)) w_i^{(1)}(D(y, F)) h(y, x) \, dF(y) \right. \right.$$
$$(17) \qquad\qquad\qquad\qquad \left. + (x - L_i(F)) w_i(D(x, F)) \right) d\nu_n(x) \right\}$$
$$\times \left\{ \int w_i(D(x, F)) \right\}^{-1} + o_p(1).$$

Note that $s_2(F_n) = s_2(F) + o_p(1)$ and $L_i(F_n) = L_i(F) + o_p(1)$, $i = 1, 2$ (see the proof of Theorem 2.1). By (12) and (15)–(17), we have

$$
\begin{aligned}
(18) \quad &\sqrt{n}((S(F_n)) - (S(F))) \\
&= \int (K_s(x, F) - K_1(x, F)(L_2(F) - L_1(F))' \\
&\qquad - (L_2(F) - L_1(F))(K_1(x, F))') \, d\nu_n(x) + o_p(1).
\end{aligned}
$$

Note that $\mathrm{vec}(ab') = b \otimes a$ for any $a, b \in \mathbb{R}^d$. The desired result now follows from the CLT. $\square$

PROOF OF THEOREM 2.3. The proof follows closely that of Theorem 2.2 and is thus omitted. $\square$

PROOF OF COROLLARY 3.1. By Theorem 3.2, $K(x) = K_s(x, F)$ since $L_i(F) = \theta$ for $i = 1, 2$. Assume without loss of generality that $\theta = 0$. For the given $F$ and $(\mu, \sigma)$, it follows that

$$
u(x) = \Sigma^{-1} x / \|\Sigma^{-1} x\|, \qquad (x \neq 0),
$$

$$
\sigma(F_{u(x)}) = m_0 \sqrt{u(x)' \Sigma u(x)}, \qquad O(x, F) = \|\Sigma^{-1/2} x\| / m_0.
$$

Let $u = z / \|z\|$. Observe that

$$
\begin{aligned}
PWS(F) &= \frac{\int x x' w_2(PD(x, F)) \, dF}{\int w_2(PD(x, F)) \, dF} = \frac{\Sigma^{1/2} (\int z z' w_2(s_0(\|z\|)) \, dF_0) \Sigma^{1/2}}{\int w_2(s_0(\|z\|)) \, dF_0} \\
&= E(\|Z\|^2 w_2(s_0(\|Z\|))) \Sigma^{1/2} \left( \int u' u \, dF_0 \right) \Sigma^{1/2} / c_0 = c_1 \Sigma
\end{aligned}
$$

by, for example, Lemma 5.1 of Lopuhaä (1989). By Lemma 3.4, it follows that for any $x, y \in \mathbb{R}^d$,

$$
f_1(x, u(y)) = \frac{\sqrt{u(y)' \Sigma u(y)}}{2p(0)} \mathrm{sign}(y' \Sigma^{-1} x),
$$

$$
f_2(x, u(y)) = \frac{\sqrt{u(y)' \Sigma u(y)}}{4p(m_0)} \mathrm{sign}(|y' \Sigma^{-1} x| - m_0 \|\Sigma^{-1/2} y\|).
$$

Note that $f_1(x, u(y))$ is an odd function of $y$. By Lemma 3.2, we have

$$
\begin{aligned}
c_0 K_s&(x, F) \\
&= \int (yy' - c_1 \Sigma) w_2^{(1)}(s_0(\|\Sigma^{-1/2} y\|)) h(y, x) \, dF(y) \\
&\quad + (xx' - c_1 \Sigma) w_2(s_0(\|\Sigma^{-1/2} y\|)) \\
&= \int \frac{(yy' - c_1 \Sigma) w_2^{(1)}(s_0(\|\Sigma^{-1/2} y\|)) O(y, F) f_2(x, u(y))}{\sigma(F_{u(y)})(1 + O(y, F))^2} \, dF(y) \\
&\quad + (xx' - c_1 \Sigma) w_2(s_0(\|\Sigma^{-1/2} x\|)).
\end{aligned}
$$

Let $\tilde{x} = \Sigma^{-1/2}x$, $\tilde{y} = \Sigma^{-1/2}y$, $\tilde{y}/\|\tilde{y}\| = u = (u_1, \ldots, u_d)'$ and $T$ be an orthogonal matrix with $\tilde{x}/\|\tilde{x}\|$ as its first column. We have

$$c_0 K_s(x, F) - (xx' - c_1\Sigma)w_2(s_0(\|\Sigma^{-1/2}x\|))$$

$$= \int \frac{(yy' - c_1\Sigma)w_2^{(1)}(s_0(\|\tilde{y}\|))\|\tilde{y}\|s_0^2(\|\tilde{y}\|)\operatorname{sign}(|(\tilde{y})'\tilde{x}| - m_0\|\tilde{y}\|)}{4m_0^2 p(m_0)}\, dF(y)$$

$$= \Sigma^{1/2}\int \{(\tilde{y}\tilde{y}' - c_1 I_d)w_2^{(1)}(s_0(\|\tilde{y}\|))\|\tilde{y}\|s_0^2(\|\tilde{y}\|)\operatorname{sign}(|(\tilde{y})'\tilde{x}| - m_0\|\tilde{y}\|)\}$$

$$\qquad\qquad \times \{4m_0^2 p(m_0)\}^{-1}\, dF_0(\tilde{y})\Sigma^{1/2}$$

$$= \Sigma^{1/2}T\int\{(\tilde{y}/\|\tilde{y}\|\tilde{y}'/\|\tilde{y}\|\|\tilde{y}\|^2 - c_1 I_d)$$

$$\qquad\qquad \times w_2^{(1)}(s_0(\|\tilde{y}\|))\|\tilde{y}\|s_0^2(\|\tilde{y}\|)\operatorname{sign}(|u_1|\|\tilde{x}\| - m_0)\}$$

$$\qquad\qquad \times \{4m_0^2 p(m_0)\}^{-1}\, dF_0(\tilde{y})T'\Sigma^{1/2}$$

$$= \Sigma^{1/2}T\Big(c_3\int uu'\operatorname{sign}(|u_1|\|\tilde{x}\| - m_0)\, dF_0(\tilde{y}) - c_1 c_2 s_1(\|\tilde{x}\|)I_d\Big)T'\Sigma^{1/2},$$

by Theorem 1.5.6 of Murihead (1982). Note that

$$Tc_3\int uu'\operatorname{sign}(|u_1|\|\tilde{x}\| - m_0)\, dF_0(\tilde{y})T'$$

$$= Tc_3\operatorname{diag}(s_2(\|\tilde{x}\|), \tilde{s}_2(\|\tilde{x}\|), \ldots, \tilde{s}_2(\|\tilde{x}\|))T'$$

$$= c_3\tilde{s}_2(\|\tilde{x}\|)I_d + c_3(s_2(\|\tilde{x}\|) - \tilde{s}_2(\|\tilde{x}\|))\frac{\tilde{x}}{\|\tilde{x}\|}\frac{\tilde{x}'}{\|\tilde{x}\|},$$

where $\tilde{s}_2(t) = \int u_2^2\operatorname{sign}(|u_1|t - m_0)\, dF_0(\tilde{y}) = (s_1(t) - s_2(t))/(d-1)$. Therefore, we have

$$K(X) = K_s(X, F) = \frac{1}{c_0}\Sigma^{1/2}\Big(t_1(\|\tilde{X}\|)\frac{\tilde{X}}{\|\tilde{X}\|}\frac{\tilde{X}'}{\|\tilde{X}\|} + t_2(\|\tilde{X}\|)I_d\Big)\Sigma^{1/2}.$$

Now invoking Lemmas 5.1 and 5.2 of Lopuhaä (1989), we obtain the desired result. $\square$

PROOF OF THEOREM 3.5. We need the following lemma. Its proof is skipped. Note that $F(\varepsilon, \delta_y) = (1-\varepsilon)F + \varepsilon\delta_y$ and $F_u(\varepsilon, \delta_y) = (1-\varepsilon)F_u + \varepsilon\delta_{u'y}$ for any unit vector $u$.

LEMMA 5.1.  *Suppose that $X \sim F$ is elliptically symmetric about the origin with a positive definite matrix $\Sigma$ associated. Let $a(u) = \sqrt{u'\Sigma u}$. Then:*

1. $\operatorname{Med}(F_u(\varepsilon, \delta_x)) = \operatorname{Med}\{-a(u)d_1(\varepsilon), u'x, a(u)d_1(\varepsilon)\}$, *and*

2. $\mathrm{MAD}(F_u(\varepsilon,\delta_x)) = \mathrm{Med}\{a(u)m_1(\mathrm{Med}(F_u(\varepsilon,\delta_x))/a(u),\varepsilon), |u'x - \mathrm{Med}(F_u(\varepsilon,\delta_x))|, a(u)m_2(\mathrm{Med}(F_u(\varepsilon,\delta_x))$

We now turn to the proof of Theorem 3.5. By Lemma 5.1, for any $y \in \mathbb{R}^d$, we have that

$$\mu(F_u(\varepsilon,\delta_y))/a(u) = \mathrm{Med}\{-a(u)d_1, u'y, a(u)d_1\}/a(u)$$
$$= \mathrm{Med}\{-d_1(\varepsilon), (\Sigma^{1/2}u)'/a(u)\Sigma^{-1/2}y, d_1(\varepsilon)\},$$
$$\frac{\sigma(F_u(\varepsilon,\delta_y))}{a(u)} = \mathrm{Med}\left\{m_1\left(\frac{\mu(F_u(\varepsilon,\delta_y))}{a(u)},\varepsilon\right),\right.$$
$$\left|\frac{(\Sigma^{1/2}u)'}{a(u)}\Sigma^{-1/2}y - \frac{\mu(F_u(\varepsilon,\delta_y))}{a(u)}\right|,$$
$$\left. m_2\left(\frac{\mu(F_u(\varepsilon,\delta_y))}{a(u)},\varepsilon\right)\right\}.$$

Let $v = \Sigma^{1/2}u/a(u)$, $\tilde{y} = \Sigma^{-1/2}y$ and $\tilde{x} = \Sigma^{-1/2}x$. Then all the mappings are one-to-one and $\|v\| = 1$. Denote $f_5(u,x,d_1) = \mathrm{Med}\{-d_1, u'x, d_1\}$. Then

$$O(x, F(\varepsilon,\delta_y))$$
$$= \sup_{\|v\|=1} \frac{v'\tilde{x} - f_5(v,\tilde{y},d_1)}{\mathrm{Med}\{m_1(f_5(v,\tilde{y},d_1),\varepsilon), |v'\tilde{y} - f_5(v,\tilde{y},d_1)|, m_2(f_5(v,\tilde{y},d_1),\varepsilon)\}}.$$

Let $U$ be an orthogonal matrix with $\tilde{y}/\|\tilde{y}\|$ as its first column, and $U'v = \tilde{v}$. Then $f_5(v,\tilde{y},d_1) = \mathrm{Med}\{-d_1, \tilde{v}_1\|\tilde{y}\|, d_1\} = f_4(\tilde{v}_1, \|\tilde{y}\|, d_1)$ and $O(x, F(\varepsilon,\delta_y))$ becomes

$$\sup_{\|\tilde{v}\|=1} \{\tilde{v}'U'\tilde{x} - f_4(\tilde{v}_1, \|\tilde{y}\|, d_1)\}$$
$$\times \{\mathrm{Med}\{m_1(f_4(\tilde{v}_1, \|\tilde{y}\|, d_1),\varepsilon),$$
$$|\tilde{v}_1\|\tilde{y}\| - f_4(\tilde{v}_1, \|\tilde{y}\|, d_1)|, m_2(f_4(\tilde{v}_1, \|\tilde{y}\|, d_1),\varepsilon)\}\}^{-1}$$
$$= \sup_{\|\tilde{v}\|=1} (\tilde{v}'U'\tilde{x} - f_4(\tilde{v}_1, \|\tilde{y}\|, d_1))/f_3(\tilde{v}_1, \|\tilde{y}\|, d_1).$$

It follows that

$$\int xx' w_2(PD(x, F(\varepsilon,\delta_y)))\, dF(x)$$
$$= \int \Sigma^{1/2}\tilde{x}\tilde{x}' w_2\left(1\Big/\left(1 + \sup_{\|\tilde{v}\|=1} \frac{\tilde{v}'U'\tilde{x} - f_4(\tilde{v}_1, \|\tilde{y}\|, d_1)}{f_3(\tilde{v}_1, \|\tilde{y}\|, d_1)}\right)\right)\Sigma^{1/2}\, dF(x)$$
$$= \int \Sigma^{1/2}Uxx' w_2\left(1\Big/\left(1 + \sup_{\|u\|=1} \frac{u'x - f_4(u_1, \|\tilde{y}\|, d_1)}{f_3(u_1, \|\tilde{y}\|, d_1)}\right)\right)U'\Sigma^{1/2}\, dF_0(x).$$

Observe that

$$\sup_{\|u\|=1} \frac{u'x - f_4(u_1, \|\tilde{y}\|, d_1)}{f_3(u_1, \|\tilde{y}\|, d_1)}$$

$$= \sup_{-1 \leq u_1 \leq 1} \sup_{\|u_2\| = \sqrt{1-u_1^2}} \frac{u_2'x_2 + u_1 x_1 - f_4(u_1, \|\tilde{y}\|, d_1)}{f_3(u_1, \|\tilde{y}\|, d_1)}$$

$$= \sup_{0 \leq u_1 \leq 1} \frac{\sqrt{1-u_1^2}\|x_2\| + |u_1 x_1 - f_4(u_1, \|\tilde{y}\|, d_1)|}{f_3(u_1, \|\tilde{y}\|, d_1)}$$

$$= f_1(x, \|\tilde{y}\|, \varepsilon),$$

which is an even function of $x_2$. Hence,

$$\int xx' w_2(PD(x, F(\varepsilon, \delta_y)))\, dF(x)$$

$$= \int \Sigma^{1/2} U xx' w_2(1/(1 + f_1(x, \|\tilde{y}\|, \varepsilon))) U' \Sigma^{1/2}\, dF_0(x)$$

$$= \Sigma^{1/2}\Big(\psi_2(\|\tilde{y}\|, \varepsilon) I_d + (\psi_1(\|\tilde{y}\|, \varepsilon) - \psi_2(\|\tilde{y}\|, \varepsilon)) \frac{\tilde{y}}{\|\tilde{y}\|} \frac{(\tilde{y})'}{\|\tilde{y}\|}\Big) \Sigma^{1/2}.$$

Likewise, we can show that

$$\int x w_i(PD(x, F(\varepsilon, \delta_y)))\, dF(x)$$

$$= (y/\|\tilde{y}\|) \int x_1 w_i(1/(1 + f_1(x, \|\tilde{y}\|, \varepsilon)))\, dF_0(x).$$

Thus

$$L_i(F(\varepsilon, \delta_y))$$

$$= \Big\{ (y/\|\tilde{y}\|)\Big( (1 - \varepsilon) \int x_1 w_i(1/(1 + f_1(x, \|\tilde{y}\|, \varepsilon)))\, dF_0(x)$$

$$+ \varepsilon\|\tilde{y}\| w_i(1/(1 + f_2(\|\tilde{y}\|, \varepsilon))) \Big) \Big\}$$

$$\times \Big\{ (1 - \varepsilon) \int w_i(1/(1 + f_1(x, \|\tilde{y}\|, \varepsilon)))\, dF_0(x)$$

$$+ \varepsilon w_i(1/(1 + f_2(\|\tilde{y}\|, \varepsilon))) \Big\}^{-1}$$

and

$$PWS(F(\varepsilon, \delta_y))$$

$$= \Big\{ (1 - \varepsilon)\Big( \psi_2(\|\tilde{y}\|, \varepsilon)\Sigma + (\psi_1(\|\tilde{y}\|, \varepsilon) - \psi_2(\|\tilde{y}\|, \varepsilon)) \frac{y}{\|\tilde{y}\|} \frac{y'}{\|\tilde{y}\|} \Big)$$

$$+ \varepsilon \|\tilde{y}\|^2 w_2 \Big( \frac{1}{1 + f_2(\|\tilde{y}\|, \varepsilon)} \Big) \frac{y}{\|\tilde{y}\|} \frac{y'}{\|\tilde{y}\|} \Big\}$$

$$\times \Big\{ (1 - \varepsilon) \int w_2 \Big( \frac{1}{1 + f_1(x, \|\tilde{y}\|, \varepsilon)} \Big) dF_0(x)$$

$$+ \varepsilon w_2 \Big( \frac{1}{1 + f_2(\|\tilde{y}\|, \varepsilon)} \Big) \Big\}^{-1}$$

$$- L_1(F(\varepsilon, \delta_y))(L_2(F(\varepsilon, \delta_y)))' - L_2(F(\varepsilon, \delta_y))(L_1(F(\varepsilon, \delta_y)))'$$

$$+ L_1(F(\varepsilon, \delta_y))(L_1(F(\varepsilon, \delta_y)))'.$$

The desired result follows. $\square$

PROOF OF THEOREM 3.6. (a) is trivial. We now show (b). Assume, w.l.o.g. that $\theta = 0$. Since $\mathrm{CSI}(PWS, F) \geq \mathrm{GESI}(PWS, F)$, we need to show that $\mathrm{CSI}(PWS, F) \leq \mathrm{GESI}(PWS, F)$. Following the proof of Theorem 2.3 and noting that $L_i(F(\varepsilon, \delta_y)) = L_i(F) + o(1)$, $i = 1, 2$, we can show that

$$(PWS(F(\varepsilon, \delta_y)) - PWS(F))/\varepsilon$$

$$= \Big( \int x x' w_2^{(1)}(PD(x, F)) H_\varepsilon(x, y) F(dx) + y y' w_2(PD(y, F)) \Big)$$

$$\times \Big( \int w_2(PD(x, F)) F(dx) \Big)^{-1}$$

$$- PWS(F) \Big( \int w_2^{(1)}(PD(x, F)) H_\varepsilon(x, y) F(dx) + w_2(PD(y, F)) \Big)$$

$$\times \Big( \int w_2(PD(x, F)) F(dx) \Big)^{-1} + o(1),$$

where $o(\cdot)$ is in the uniform sense with respect to $y \in \mathbb{R}^d$. Following the proof of Theorem 3.5 of Zuo, Cui and Young (2004) and letting $g(x, u, F) = (u, x - \mu(F_u))/\sigma(F_u)$, we have

$$\Big( \int x x' w_2^{(1)}(PD(x, F)) H_\varepsilon(x, y) F(dx) + y y' w_2(PD(y, F)) \Big)$$

$$\times \Big( \int w_2(PD(x, F)) F(dx) \Big)^{-1}$$

$$= \Big\{ \int_{S(x, M)} x x' w_2^{(1)}(PD(x, F)) \frac{g(x, u(x), F) - g(x, u(x), F(\varepsilon, \delta_y))}{\varepsilon(1 + O(x, F))^2} dF(x)$$

$$+ y y' w_2(PD(y, F)) \Big\}$$

$$\times \Big\{ \int w_2(PD(x, F)) dF(x) \Big\}^{-1} + I_5(M, y, \varepsilon) + o(1),$$

where $S(x, M) = \{x : 1/M \leq \|\Sigma^{-1/2}x\| \leq M\}$ for a fixed $M > 0$, $\sup_{y \in \mathbb{R}^d, \varepsilon < 0.5} \|I_5(M, y, \varepsilon)\| \to 0$ as $M \to \infty$ and $o(\cdot)$ is in the uniform sense in $y \in \mathbb{R}^d$. Note that $u(x) = \Sigma^{-1}x/\|\Sigma^{-1}x\|$ and $\sigma(F_{u(x)}) = m_0\|\Sigma^{-1/2}x\|/\|\Sigma^{-1}x\|$ for $x \neq 0$, $\mu(F_{u(x)}) = 0$ and $O(x, F) = \|\Sigma^{-1/2}x\|/m_0$. By Lemma 5.1 we see that $\mu(F(\varepsilon, \delta_y))$ is odd in $y$. Therefore,

$$
\left( \int xx'w_2^{(1)}(D(x, F))H_\varepsilon(x, y)F(dx) + yy'w_2(D(y, F)) \right)
$$
$$
\times \left( \int w_2(D(x, F))F(dx) \right)^{-1}
$$
$$
= \left\{ \int_{S(x,M)} xx'w_2^{(1)}(s_0(\|\Sigma^{-1/2}x\|))\|\Sigma^{-1}x\| \right.
$$
$$
\left. \times \frac{\sigma(F_{u(x)}(\varepsilon, \delta_y)) - \sigma(F_{u(x)})}{m_0^2\varepsilon(1 + O(x, F))^2} dF(x) + yy'w_2(PD(y, F)) \right\}
$$
$$
\times \left\{ \int w_2(PD(x, F)) dF(x) \right\}^{-1} + I_5(M, y, \varepsilon) + o(1),
$$

where $o(\cdot)$ is in the uniform sense with respect to $y \in \mathbb{R}^d$. Call the first term in the right-hand side of the last equality $I_6 = I_6(M, y, \varepsilon)$. By Lemma 5.1,

$$
\sigma(F_{u(x)}(\varepsilon, \delta_y)) - \sigma(F_{u(x)})
$$
$$
= \frac{\|\Sigma^{-1/2}x\|}{\|\Sigma^{-1}x\|} \mathrm{Med}\left\{ m_1\left( \frac{\mu(F(\varepsilon, \delta_y))}{a(u(x))}, \varepsilon \right), \right.
$$
$$
\left. \left| \frac{x'\Sigma^{-1}y}{\|\Sigma^{-1/2}x\|} - \frac{\mu(F(\varepsilon, \delta_y))}{a(u(x))} \right|, m_2\left( \frac{\mu(F(\varepsilon, \delta_y))}{a(u(x))}, \varepsilon \right) \right\},
$$

where $\mu(F(\varepsilon, \delta_y))/a(u(x)) = \mathrm{Med}\{-d_1, x'\Sigma^{-1}y/\|\Sigma^{-1/2}x\|, d_1\}$. Let $\tilde{x} = \Sigma^{-1/2}x$, $\tilde{y} = \Sigma^{-1/2}y$ and $T$ be an orthogonal matrix with $\tilde{y}/\|\tilde{y}\|$ as its first column. Note that $T'\tilde{X} \overset{d}{=} \tilde{X}$. Denote $T'\tilde{x}/\|\tilde{x}\| = u = (u_1, \ldots, u_d)'$. Changing variables $(\tilde{x} = \Sigma^{-1/2}x)$ and then taking an orthogonal transformation (with matrix $T$) and taking advantage of the independence of $\|\tilde{X}\|$ and $\tilde{X}/\|\tilde{X}\|$ [see Lemma 5.1 of Lopuhaä (1999)], we have

$$
I_6 = \left\{ \Sigma^{1/2}T \int_{1/M \leq \|\tilde{x}\| \leq M} \tilde{x}\tilde{x}'w_2^{(1)}(s_0(\|\tilde{x}\|))\|\tilde{x}\|s_0^2(\|\tilde{x}\|) \frac{I_7(u_1, \tilde{y}, \varepsilon)}{m_0^2\varepsilon} dF_0(\tilde{x}) \right.
$$
$$
\left. \times T'\Sigma^{1/2} + yy'w_2(s_0(\|\tilde{y}\|)) \right\}
$$
$$
\times \left\{ \int w_2(s_0(\|\tilde{x}\|)) dF_0(\tilde{x}) \right\}^{-1}
$$

$$\begin{aligned}
&= \Big\{ \Sigma^{1/2} T \int_{1/M \leq \|\tilde{x}\| \leq M} \|\tilde{x}\|^3 w_2^{(1)}(s_0(\|\tilde{x}\|)) s_0^2(\|\tilde{x}\|) \, dF_0(\tilde{x}) \\
&\qquad\qquad \times \int_{1/M \leq \|\tilde{x}\| \leq M} uu' \frac{I_7(u_1, \tilde{y}, \varepsilon)}{m_0^2 \varepsilon} \, dF_0(\tilde{x}) T' \Sigma^{1/2} \Big\} \\
&\quad \times \Big\{ \int w_2(s_0(\|\tilde{x}\|)) \, dF_0(\tilde{x}) \Big\}^{-1} \\
&\quad + yy' w_2(s_0(\|\tilde{y}\|)) \Big/ \int w_2(s_0(\|\tilde{x}\|)) \, dF_0(\tilde{x}),
\end{aligned}$$

where $I_7(u_1, \tilde{y}, \varepsilon) = \mathrm{Med}\{m_1(I_8(u_1, \tilde{y}, \varepsilon), \varepsilon), |u_1\|\tilde{y}\| - I_8(u_1, \tilde{y}, \varepsilon)|, m_2(I_8(u_1, \tilde{y}, \varepsilon), \varepsilon)\} - m_0$ and $I_8(u_1, \tilde{y}, \varepsilon) = \mathrm{Med}\{-d_1, u_1\|\tilde{y}\|, d_1\}$. It can be shown (details are skipped) that

$$I_7(u_1, \tilde{y}, \varepsilon)/\varepsilon = \mathrm{sign}(|u_1\|\tilde{y}\| - m_0|)/(4p(m_0)) + o(1),$$

where $o(1) \to 0$ uniformly in $y \in \mathbb{R}^d$ as $\varepsilon \to 0$. Following the proof of Corollary 3.1, we have

$$\begin{aligned}
I_6 &= \frac{\Sigma^{1/2} T c_3(M) \int_{1/M \leq \|\tilde{x}\| \leq M} uu' \, \mathrm{sign}(|u_1\|\tilde{y}\| - m_0|) T' \Sigma^{1/2}}{\int w_2(s_0(\|\tilde{x}\|)) \, dF_0(\tilde{x})} \\
&\quad + \frac{yy' w_2(s_0(\|\tilde{y}\|))}{\int w_2(s_0(\|\tilde{x}\|)) \, dF_0(\tilde{x})} + o(1) \\
&= \frac{\Sigma^{1/2} c_3(M)(\tilde{s}_2(\|\tilde{y}\|, M) I_d + (s_2(\|\tilde{y}\|, M) - \tilde{s}_2(\|\tilde{y}\|, M)) \tilde{y}/\|\tilde{y}\| \tilde{y}'/\|\tilde{y}\|) \Sigma^{1/2}}{\int w_2(s_0(\|\tilde{x}\|)) \, dF_0(\tilde{x})} \\
&\quad + yy' w_2(s_0(\|\tilde{y}\|)) \Big/ \int w_2(s_0(\|\tilde{x}\|)) \, dF_0(\tilde{x}) + o(1),
\end{aligned}$$

where $o(1)$ is in the same sense as before. Further,

$$c_3(M) = \int_{1/M \leq \|\tilde{x}\| \leq M} \|\tilde{x}\|^3 w_2^{(1)}(s_0(\|\tilde{x}\|)) s_0^2(\|\tilde{x}\|) \, dF_0(\tilde{x}),$$

$$s_2(t, M) = \int_{1/M \leq \|\tilde{x}\| \leq M} u_1^2 \, \mathrm{sign}(|u_1 t - m_0|) \, dF_0(\tilde{x}),$$

$$\tilde{s}_2(t, M) = \int_{1/M \leq \|\tilde{x}\| \leq M} u_2^2 \, \mathrm{sign}(|u_1 t - m_0|) \, dF_0(\tilde{x}).$$

Therefore,

$$\begin{aligned}
&(PWS(F(\varepsilon, \delta_y)) - PWS(F))/\varepsilon \\
&\quad = \Sigma^{1/2} \frac{\tilde{y}}{\|\tilde{y}\|} (c_3(M)(s_2(\|\tilde{y}\|, M) - \tilde{s}_2(\|\tilde{y}\|, M)) + \|\tilde{y}\|^2 w_2(s_0(\|\tilde{y}\|))) \frac{\tilde{y}'}{\|\tilde{y}\|} \Sigma^{1/2}
\end{aligned}$$

$$\times \left( \int w_2(s_0(\|\tilde{x}\|)) \, dF_0(\tilde{x}) \right)^{-1}$$

$$+ I_5(M, y, \varepsilon) + o(1)$$

$$+ \Sigma \left( \frac{c_3(M)\tilde{s}_2(\|\tilde{y}\|, M)}{\int w_2(s_0(\|\tilde{x}\|)) \, dF_0(\tilde{x})} \right.$$

$$\left. - c_1 \frac{\int w_2^{(1)}(D(x,F)) H_\varepsilon(x,y) F(dx) + w_2(D(y,F))}{\int w_2(D(x,F)) F(dx)} \right),$$

where again $o(1) \to 0$ uniformly in $y \in \mathbb{R}^d$ as $\varepsilon \to 0$. From the definition of CSI, it follows that

$$\text{CSI}(PMS, F)$$

$$= \lim_{\varepsilon \to 0^+} \sup_{y \in \mathbb{R}^d} \left\| \left\{ \Sigma^{1/2} \tilde{y}/\|\tilde{y}\| (c_3(M)(s_2(\|\tilde{y}\|, M) - \tilde{s}_2(\|\tilde{y}\|, M)) \right. \right.$$

$$\left. + \|\tilde{y}\|^2 w_2(s_0(\|\tilde{y}\|)) \tilde{y}'/\|\tilde{y}\| \Sigma^{1/2} \right\}$$

$$\left. \times \left\{ \int w_2(s_0(\|\tilde{x}\|)) \, dF_0(\tilde{x}) \right\}^{-1} + I_5(M, y, \varepsilon) + o(1) \right\|$$

$$\leq \lim_{\varepsilon \to 0^+} \sup_{y \in \mathbb{R}^d} \left\| \left\{ \Sigma^{1/2} \tilde{y}/\|\tilde{y}\| (c_3(M)(s_2(\|\tilde{y}\|, M) - \tilde{s}_2(\|\tilde{y}\|, M)) \right. \right.$$

$$\left. + \|\tilde{y}\|^2 w_2(s_0(\|\tilde{y}\|)) \tilde{y}'/\|\tilde{y}\| \Sigma^{1/2} \right\}$$

$$\left. \times \left\{ \int w_2(s_0(\|\tilde{x}\|)) \, dF_0(\tilde{x}) \right\}^{-1} \right\|$$

$$+ \lim_{\varepsilon \to 0^+} \sup_{y \in \mathbb{R}^d} \| I_5(M, y, \varepsilon) \|$$

$$\leq \lambda_1 \sup_{r \geq 0} \left| \frac{c_3(M)(s_2(r, M) - \tilde{s}_2(r, M)) + r^2 w_2(s_0(r))}{\int w_2(s_0(\|\tilde{x}\|)) \, dF_0(\tilde{x})} \right|$$

$$+ \lim_{\varepsilon \to 0^+} \sup_{y \in \mathbb{R}^d} \| I_5(M, y, \varepsilon) \|.$$

Now letting $M \to \infty$, we get $\text{CSI}(PMS, F) \leq \lambda_1 \sup_{r \geq 0} |t_1(r)|/c_0 = \text{GESI}(PMS, F)$. $\square$

**Acknowledgments.** We appreciate the insightful remarks and the constructive criticism of an Associate Editor, Editor John Marden and two referees, which led to improvements in the paper. We also thank Xuming He and Ricardo Maronna for helpful suggestions.

Department of Statistics and Probability
Michigan State University
East Lansing, Michigan 48824
USA
E-mail: zuo@msu.edu

Department of Mathematics
Beijing Normal University
Beijing 100875
China
E-mail: cuihj@hotmail.com